\documentclass{jglta}
\usepackage{amsmath}
\usepackage{amsthm}
\usepackage{amsfonts}
\usepackage{amssymb}
\usepackage{parsetree}
\JGLTAnumberwithin{equation}{section}

\allowdisplaybreaks
\newtheorem{theorem}{Theorem}
\newtheorem{proposition}[theorem]{Proposition}

\newtheorem{cor}[theorem]{Corollary}
\theoremstyle{definition}
\newtheorem{definition}{Definition}
 
\newtheorem{rk}[definition]{Remark}
\begin{document}
\renewcommand{\evenhead}{Marius Buliga}
\renewcommand{\oddhead}{Dilatation structures I. Fundamentals}
\thispagestyle{empty}
\FirstPageHead{*}{*}{20**}{\pageref{firstpage}--\pageref{lastpage}}

\label{firstpage}

\Name{Dilatation structures I. Fundamentals}
\Author{Marius BULIGA}
\Address{"Simion Stoilow" Institute of Mathematics of the Romanian Academy, \\
P.O. BOX 1-764, RO 014700  Bucure\c sti, Romania\\  
\smallskip
E-mail: Marius.Buliga@imar.ro}

\begin{abstract}
A dilatation structure is a concept in between a group and a differential structure. 
In this article  we study  fundamental properties of dilatation structures on 
metric spaces. This is a  part  of a series of papers which show that  such a structure allows
 to do non-commutative  analysis, in the sense of differential calculus, on a large class of metric spaces, 
some of them  fractals. We also describe a formal, universal calculus with binary decorated planar trees, which underlies any 
dilatation structure.
\end{abstract}

\tableofcontents

\section{Introduction}

The purpose of this paper is to introduce dilatation structures on metric spaces. 
A  dilatation structure is a concept in between a group 
and a differential structure. Any metric space $(X,d)$  endowed with a dilatation 
structure has an associated tangent bundle. The tangent space at a point is a 
conical group, that is the tangent space  has a group structure together with a one-parameter 
group of automorphisms. Conical groups generalize Carnot groups, i.e 
nilpotent groups endowed with a graduation. Each dilatation structure leads to a 
non-commutative differential calculus on the metric space $(X,d)$. 

There  are  several  important papers dedicated to the study of extra 
structures on a metric space which allows to do a reasonable analysis in such 
spaces, like Cheeger \cite{cheeger} or Margulis-Mostow \cite{marmos1,marmos2}.

The constructions proposed in this paper first  appeared  in connection to problems in analysis on 
sub-riemannian manifolds. Parts of this article can be seen  as a rigorous formulation of the 
considerations in the last section of Bella\"{\i}che \cite{bell}.

A dilatation structure is simply  a bundle of semigroups of (quasi-)contractions 
on the metric space $(X,d)$, satisfying a number of axioms.
The tangent bundle structure associated with a given dilatation structure 
on the metric space $(X,d)$  is obtained by a passage to the limit procedure, 
starting from an  algebraic structure which lives on the metric space. 

With the help of the dilatation structure we construct a bundle 
(over the metric space) of (local) operations: to each $x\in X$ and parameter 
$\varepsilon$, for simplicity here $\varepsilon \in (0,+\infty)$, 
there is a natural non-associative  operation 
\[\Sigma^{x}_{\varepsilon}: U(x) \times U(x) \rightarrow U(x) \]
where $U(x)$ is a neighbourhood of $x$. The  non-associativity of this 
operation is controlled by the parameter 
$\varepsilon$. As $\varepsilon$ goes to $0$ the operation 
$\displaystyle \Sigma^{x}_{\varepsilon}$ converges to a group operation 
on the tangent space of $(X,d)$ at $x$. 

 Denote by 
$\displaystyle \delta^{x}_{\varepsilon}$ the dilatation based at $x \in X$, 
of parameter $\varepsilon$. The bundle of operations  satisfies a kind of weak 
associativity, even if for any fixed $y \in X$ the operation   
$\displaystyle \Sigma^{y}_{\varepsilon}$ is non-associative. The weak associativity 
property, named also shifted associativity, is
\[\Sigma^{x}_{\varepsilon}(u, \Sigma^{\delta^{x}_{\varepsilon}}_{\varepsilon} 
(v,w)) = \Sigma^{x}_{\varepsilon}(\Sigma^{x}_{\varepsilon}(u,v), w) \]
for any $x \in X$ and any $u,v,w \in X$ sufficiently close to $X$. We shall 
describe  also other 
objects (like a function satisfying a shifted inverse property) and 
algebraic identities related to the dilatation structure and 
the induced bundle of operations.

We briefly describe further the contents of the paper. In section 2 we give 
motivational examples of dilatation structures. Basic notions and results of 
metric geometry and groups endowed with dilatations  are mentioned  
in section 3. 

In section 4 we introduce a formalism based on decorated planar binary trees. 
This formalism will be used to prove the main results of the paper. We show that, 
from an algebraic point of view, dilatation structures (more precisely the 
formalism in section 4) induce a bundle of one parameter deformations of 
 binary operations, which are not associative, but shifted associative. 
This is a structure which bears resemblance with the tangent bundle of a  Lie 
group, but it is more general. 

Section 5, 6 and 7  are devoted to dilatation structures. These sections 
contain the main results of the paper. After we introduce and explain the 
axioms of dilatation structures, we describe  several key metric properties of 
such a  structure, in section 5. In section 6 we prove that a dilatation 
structure induces a valid notion of tangent bundle. 
In section 7 we explain how a dilatation structure leads to a  differential 
calculus. 

Section 8 is made of two parts. In the first part we show that dilatation 
structures induce differential structures, in a generalized sense. In the 
second part we turn to conical groups and  we prove the curious result that, 
even if in a conical group left translations are smooth but right translations 
are generically non differentiable, the group operation 
is smooth if we well choose a dilatation structure.


\section{Motivation}

We start with a trivial example of a dilatation structure, then we 
briefly explain the occurence of such a structure in more unusual 
situations. 

There is a lot of structure hiding in the dilatations of 
$\displaystyle \mathbb{R}^{n}$. For this space, the dilatation based 
at $x$, of coefficient $\varepsilon>0$, is the function 
\[\delta^{x}_{\varepsilon}: \mathbb{R}^{n}\rightarrow \mathbb{R}^{n} \quad 
\delta^{x}_{\varepsilon} y = x + \varepsilon (-x+y) \]
For fixed $x$ the dilatations based at $x$ form a one parameter group which
 contracts any bounded neighbourhood of $x$ to a point, uniformly with respect 
to $x$. 

Dilatations behave well with respect to the euclidean distance $d$, in the following sense: 
for any $\displaystyle x,u,v\in \mathbb{R}^{n}$ and any $\varepsilon>0$ we have 
\[\frac{1}{\varepsilon} d( \delta_{\varepsilon}^{x} u , \delta^{x}_{\varepsilon} v ) = d(u,v) \]
This shows that from the metric point of view the space $\displaystyle 
(\mathbb{R}^{n}, d)$ is a metric cone, that is $\displaystyle 
(\mathbb{R}^{n}, d)$ looks the same at all scales. 

Moreover, let $\displaystyle f: \mathbb{R}^{n} \rightarrow \mathbb{R}^{n}$ be a function 
and $\displaystyle x\in \mathbb{R}^{n}$. The function $f$ is differentiable 
in $x$ if there is a linear transformation $A$ (that is a group morphism which 
commutes with dilatations based at the neutral element $0$) such that the 
limit 
\begin{equation}
\lim_{\varepsilon\rightarrow 0} \delta_{\varepsilon^{-1}}^{f(x)} f \delta^{x}_{\varepsilon} 
(v) = f(x) +  A(-x+v) 
\label{mot1}
\end{equation}
is uniform with respect to $v$ in bounded neighbourhood of $x$. Really, let us calculate
\[\delta_{\varepsilon^{-1}}^{f(x)} f \delta^{x}_{\varepsilon} (v) = 
f(x) + \frac{1}{\varepsilon} ( -f(x)+f(x+\varepsilon(-x+v))) \]
This shows that we get the usual definition of differentiability. 

The relation (\ref{mot1}) can be put in another form, using the euclidean distance: 
\[\lim_{\varepsilon\rightarrow 0} \frac{1}{\varepsilon} d( \delta_{\varepsilon}^{f(x)} T(x)(v),  f(\delta^{x}_{\varepsilon} 
v)) = 0 \]
uniformly with respect to $v$ in  bounded neighbourhood of $x$. 
Here 
\[
T(x)(v)= x+ A(-x+v)\]
In conclusion, dilatations are the fundamental 
object for doing differential calculus on $\displaystyle \mathbb{R}^{n}$.  

Even the algebraic structure of $\displaystyle \mathbb{R}^{n}$ is encoded in 
dilatations. Really, we can recover the operation of addition from dilatations. 
It goes like this: for $\displaystyle x,u,v \in \mathbb{R}^{n}$ and $\varepsilon>0$ 
define 
\[\Delta_{\varepsilon}^{x}(u,v) = \delta_{\varepsilon^{-1}}^{\delta_{\varepsilon}^{x} u}
 \delta^{x}_{\varepsilon} v,
\qquad 
\Sigma_{\varepsilon}^{x}(u,v) 
= \delta_{\varepsilon^{-1}}^{x} \delta_{\varepsilon}^{\delta_{\varepsilon}^{x} u}(v),
\qquad 
inv^{x}_{\varepsilon}(u) =  \delta_{\varepsilon^{-1}}^{\delta_{\varepsilon}^{x} u} x \]
For fixed $x,u,\varepsilon$ the functions $\displaystyle \Delta_{\varepsilon}^{x}(u,\cdot), 
\Sigma_{\varepsilon}^{x}(u,\cdot)$ are inverse one to another, but we don't insist 
on this for the moment (see Proposition \ref{psumdif}). 

What is the meaning of these functions? Let us compute
\begin{align*}
\Delta_{\varepsilon}^{x}(u,v) 
&= \delta_{\varepsilon}^{x} u  + \frac{1}{\varepsilon}
\left(-\left(\delta_{\varepsilon}^{x} u\right)+ \delta^{x}_{\varepsilon} v\right) \\
&= \left( x+\varepsilon(-x+u)\right) + \frac{1}{\varepsilon}  
\left(\varepsilon(-u+x)-x+x+\varepsilon(-x+v)\right)\\
&= x+\varepsilon(-x+u) + \frac{1}{\varepsilon}\varepsilon(-u+v)\\
&=x+\varepsilon(-x+u) + (-u+v) \\
\Sigma_{\varepsilon}^{x}(u,v) 
&= x+ \frac{1}{\varepsilon} \left(-x+ \delta_{\varepsilon}^{x} u + 
\varepsilon\left(-\left(\delta_{\varepsilon}^{x} u\right)+ v\right)\right) \\
&= x+ \frac{1}{\varepsilon} \left(\varepsilon(-x+u)+ \varepsilon\left(\varepsilon(-u+x) -x+v\right)\right)\\
&= u+ \varepsilon(-u+x) + (-x+v)
\end{align*}
In the same way we get 
\[inv^{x}_{\varepsilon}(u) = =x+\varepsilon(-x+u) + (-u+x)\]
As $\varepsilon \rightarrow 0$ we have the following limits: 
\begin{gather*}
\lim_{\varepsilon\rightarrow 0} \Delta_{\varepsilon}^{x}(u,v) = \Delta^{x}(u,v) = x+(-u+v)\\
\lim_{\varepsilon\rightarrow 0} \Sigma_{\varepsilon}^{x}(u,v) = \Sigma^{x}(u,v) = u+(-x+v)\\
\lim_{\varepsilon\rightarrow 0} inv^{x}_{\varepsilon}(u) = inv^{x}(u) = x-u+x
\end{gather*}
uniform with respect to $x,u,v$ in bounded sets. 
The function $\displaystyle  \Sigma^{x}(\cdot,\cdot)$ is a group operation, namely the addition operation 
translated such that the neutral element is $x$. Thus, for $x=0$, we recover the group operation. The function 
$\displaystyle inv^{x}(\cdot)$ is the inverse function, and $\displaystyle  \Delta^{x}(\cdot,\cdot)$ is the difference 
function.

Notice that for fixed $x, \varepsilon$ the function $\displaystyle  \Sigma^{x}_{\varepsilon}(\cdot,\cdot)$ is not a 
group operation, first of all because it is not associative. Nevertheless, this function satisfies 
a shifted associativity property, namely (see Proposition \ref{p34}) 
\[\Sigma_{\varepsilon}^{x}(\Sigma_{\varepsilon}^{x}(u,v),w) = 
\Sigma_{\varepsilon}^{x}(u, \Sigma_{\varepsilon}^{\delta^{x}_{\varepsilon}u}(v,w))\]
Also, the inverse function $\displaystyle inv^{x}_{\varepsilon}$ is not involutive, but shifted involutive (Proposition 
\ref{p33}), 
\[inv_{\varepsilon}^{\delta^{x}_{\varepsilon}u}\left( inv^{x}_{\varepsilon} u\right) = u \]

These and other properties of dilatations allow to recover the structure of the tangent bundle of 
$\displaystyle \mathbb{R}^{n}$, which is trivial in this case. 

Let us go to more elaborate examples. We may look to a riemannian manifold $M$, which is locally a deformation 
of $\displaystyle \mathbb{R}^{n}$. We can use charts for transporting (locally) the dilatation structure 
from $\displaystyle \mathbb{R}^{n}$ to the manifold. All the previously described metric and algebraic properties   
will hold in this situation, in a weaker form. For example the riemannian distance is no longer scalling invariant, but 
we still have 
\[\lim_{\varepsilon\rightarrow 0} 
\frac{1}{\varepsilon} d( \delta_{\varepsilon}^{x} u , \delta^{x}_{\varepsilon} v ) = d^{x}(u,v)  \]
uniform limit with respect tu $x,u,v$ in (small) bounded sets. Here $\displaystyle d^{x}$ is an euclidean 
distance which can be identified with the distance in the tangent space of $M$ at $x$, induced by the metric 
tensor at $x$. In the same way we can construct the algebraic structure of the  tangent space at $x$, using 
the functions $\displaystyle \Sigma_{\varepsilon}^{x} , \Delta_{\varepsilon}^{x}$. We will have a differentiability notion 
coming from the dilatations transported by the chart.

If we change charts or the riemannian metric then the dilatation structure will change too, but not very much, 
essentially because the change of charts is smooth, therefore we are still able to say what are tangent spaces and 
to describe their algebraic structure.

Let us go further with more complex examples. Consider the Heisenberg group $H(n)$. As a set $\displaystyle 
H(n) = \mathbb{R}^{2n}\times\mathbb{R}$. We shall use the following notation: an element of $H(n)$ will 
be denoted by $\tilde{x}=(x,\bar{x})$, with $\displaystyle x\in \mathbb{R}^{2n}, \bar{x}\in\mathbb{R}$. The group 
operation is 
\[\tilde{x}\tilde{y} = (x+y, \bar{x}+\bar{y}+2\omega(x,y)) \]
where $\omega$ is the canonic symplectic  2-form on $\displaystyle \mathbb{R}^{2n}$. 

The group $H(n)$ is nilpotent, in fact a 2 graded Carnot group. This means that $H(n)$ is nilpotent and that 
it admits a one-parameter group of isomorphisms 
\[\delta_{\varepsilon} (x,\bar{x}) = (\varepsilon x , \varepsilon^{2} \bar{x}) \]
These are dilatations, more precisely we can construct dilatations based at $\tilde{x}$ by the formula 
\[\delta_{\varepsilon}^{\tilde{x}} \tilde{u} = \tilde{x} \delta_{\varepsilon}\left(\tilde{x}^{-1}\tilde{u}\right) \]
We may also put a scalling invariant distance on $H(n)$, for example as follows:
\[d(\tilde{x}, \tilde{y}) = g(\tilde{x}^{-1}\tilde{y}),
\qquad 
g(\tilde{u}) = \max\left\{ \|u\|, \sqrt{\mid \bar{u}\mid} \right\} \]
We can repeat step by step the constructions explained before in this situation. There are some differences 
though. 

First of all, from the metric point of view, $(H(n),d)$ is a fractal space, in the sense that 
the Hausdorff dimension of this space is equal to $2n+2$, therefore strictly greater than the topological 
dimension, which is $2n+1$. Second, the differential of a function defined by the dilatations is not the usual 
differential, but an essentially  different one, called Pansu derivative (see \cite{pansu}). 
This is part of a very active area of research in geometric analysis (among fundamental references one may cite 
\cite{pansu,cheeger,marmos1,marmos2,gromovgr}). A spectacular application of 
Pansu derivative was to prove a Rademacher theorem which in turn implies deep results about Mostow rigidity. 
The theory applies to general Carnot groups. 

The Heisenberg group is not commutative. It is in fact the model for the tangent space of a contact metric 
manifold, as the euclidean $\displaystyle \mathbb{R}^{n}$ is the model of the tangent space of a riemannian manifold. 
We enter here in the realm of sub-riemannian geometry (see for example \cite{bell,gromovsr}). 
In a future paper we shall deal with dilatation structures for sub-riemannian manifolds. An important 
problem in sub-riemannian geometry is to have good tangent bundle structures, which in turn allow us 
to prove basic theorems, like Poincar\'e inequality, Rademacher or Stepanov theorems. 

We may even go further and find dilatation structures related with rectifiable sets, or with some self-similar sets. 
This is not the purpose of this paper though. In the sequel we shall define and study fundamental properties 
of dilatation structures.

\section{Basic notions}

We denote by $f \subset X \times Z$ 
a relation and we write $f(x) = y$ if $(x,y) \in f$. Therefore we may have $f(x) = y$ and $f(x) = y'$ with 
$y \not = y'$, if $(x,y) \in f$ and $(x,y') \in f$. 

The domain of $f$ is the set of $x \in X$ such that there is $z \in Z$ with $f(x) = z$. We denote the domain by $dom \ f$. The image of $f$ is the set of $z \in Z$ such that there is $x \in X$ with $f(x) =  z$. We denote the image by $im \ f$. 

By convention, when we state that a relation $R(f(x), f(y), ...)$ is true, it means that $R(x',y', ...)$ is 
true for any choice of $x', y', ...$, such that $(x,x'), (y,y'), ... \in f$. 

In a metric space $(X,d)$, the ball centered at $x \in X$ and radius $r>0$ is denoted by $B(x,r)$.  If we need to emphasize the 
dependence on the distance $d$ then we shall use the notation  $\displaystyle B_{d}(x,r)$. In the same way, $\bar{B}(x,r)$ and  
$\displaystyle \bar{B}_{d}(x,r)$ denote the closed ball centered at $x$, with radius $r$. 

We shall use  the following convenient notation: by $\mathcal{O}(\varepsilon)$ we mean a positive function such that $\displaystyle \lim_{\varepsilon \rightarrow 0} \mathcal{O}(\varepsilon) = 0$. 

 \subsection{Gromov-Hausdorff distance}

There are several definitions of distances between metric spaces. 
For this subject see  \cite{burago} (Section 7.4),  \cite{gromov} (Chapter 3) and \cite{gromovgr}.

We explain now  a well-known alternative definition of the  Gromov-Hausdorff distance, up to a multiplicative factor.  

\begin{definition}
Let $\displaystyle (X_{i},  d_{i}, x_{i})$, $i=1,2$, be a pair of locally compact pointed metric spaces and $\mu > 0$. 
We shall say that $\mu$ is admissible if  there is a relation $\displaystyle \rho \subset X_{1} \times X_{2}$ such that 
\begin{enumerate}
\itemsep-3pt
\item[1.] $dom \ \rho$ is $\mu$-dense in $\displaystyle X_{1}$, 
\item[2.] $im \ \rho$ is $\mu$-dense in $\displaystyle X_{2}$, 
\item[3.] $\displaystyle (x_{1}, x_{2}) \in \rho$, 
\item[4.] for all $x,y \in \ dom \ \rho$ we have 
\begin{equation}
\mid d_{2}(\rho(x), \rho(y)) - d_{1}(x,y) \mid \ \leq \  \mu
\label{photoquasi}
\end{equation}
\end{enumerate}
The Gromov-Hausdorff distance  between $\displaystyle (X_{1}, x_{1}, d_{1})$ and $\displaystyle  (X_{2}, x_{2}, d_{2})$  is   the infimum of admissible numbers $\mu$. 
\label{defgh}
\end{definition}

Denote by $[X,d_{X}, x]$ the isometry class of $(X,d_{X}, x)$, that is the class of spaces $(Y,d_{Y}, y)$ such 
that it exists an isometry  $f:X \rightarrow Y$ with the property $f(x)=y$. Note that if  $(X,d_{X}, x)$ is 
isometric with $(Y,d_{Y}, y)$ then they have the same diameter. 

The Gromov-Hausdorff distance is in fact almost a  distance between isometry classes of pointed metric spaces. 
Indeed, if two pointed metric spaces are isometric then the Gromov-Hausdorff distance  equals $0$. 
The converse is also true  in the class of compact (pointed) metric spaces  \cite{gromov} (Proposition 3.6).

Moreover, if two  of the isometry classes  $[X,d_{X}, x]$, $[Y,d_{Y}, y]$, $[Z,d_{Z}, z]$ have 
(representants with) diameter at most equal to 3, then the triangle inequality is true. We  shall 
use this distance and the induced  convergence for isometry classes of the form $[X,d_{X}, x]$, 
with $diam \ X \ \leq 5/2$. 

\subsection{Metric profiles. Metric tangent space}

We shall denote by $CMS$ the set of isometry classes of pointed compact metric spaces. The distance on this set is the Gromov distance between (isometry classes of) pointed metric spaces and the topology 
is induced by this distance. 

To any locally compact metric  space we can associate a metric profile \cite{buliga3,buliga4}. 

\begin{definition}
\label{dmprof}
The metric profile associated to the locally metric space $(M,d)$ is  the assignment 
(for small enough $\varepsilon > 0$) 
\[(\varepsilon > 0 , \ x \in M) \ \mapsto \  \mathbb{P}^{m}(\varepsilon, x) = \left[\bar{B}(x,1), 
\frac{1}{\varepsilon} d, x\right] \in CMS\]
\end{definition}

We can define a notion of metric profile regardless to any distance. 

\begin{definition}
\label{dprofile}
A metric profile is a curve $\mathbb{P}:[0,a] \rightarrow CMS$ such that
\begin{enumerate}
\itemsep-3pt
\item[(a)] it is continuous at $0$,
\item[(b)]for any $b \in [0,a]$ and  $\varepsilon \in (0,1]$ we have 
\[d_{GH} (\mathbb{P}(\varepsilon b), \mathbb{P}^{m}_{d_{b}}(\varepsilon,x_{b}))  \  = \ O(\varepsilon)\]
\end{enumerate}
The function $\mathcal{O}(\varepsilon)$ may change with  $b$.
We used the notations
\[
\mathbb{P}(b) = [\bar{B}(x,1) ,d_{b}, x_{b}] \quad \mbox{  and } \quad 
\mathbb{P}^{m}_{d_{b}}(\varepsilon,x) = \left[\bar{B}(x,1),\frac{1}{\varepsilon}d_{b}, x_{b}\right] 
\]
The metric profile is nice if 
\[d_{GH}\left (\mathbb{P}(\varepsilon b), \mathbb{P}^{m}_{d_{b}}(\varepsilon,x)\right) =O(b \varepsilon) \]
\end{definition}

Imagine that $1/b$ represents the magnification on the scale of a microscope. We use the microscope to study a specimen. For each $b > 0$ the information that we get is the table of distances of the pointed metric space $\displaystyle (\bar{B}(x,1) ,d_{b}, x_{b})$. 

How can we know, just from the information given by the microscope, that the string of "images" that we have corresponds to a real specimen? The answer is that a reasonable check is the relation from point (b) of the definition of metric profiles \ref{dprofile}. 

Really, this point says that starting from any magnification $1/b$,  if we further select the ball 
$\displaystyle \bar{B}(x, \varepsilon)$ in the snapshot  $\displaystyle (\bar{B}(x,1) ,d_{b}, x_{b})$, then 
the metric space $\displaystyle (\bar{B}(x,1) ,\frac{1}{\varepsilon} d_{b}, x_{b})$ looks approximately the same as the snapshot $\displaystyle (\bar{B}(x,1) ,d_{b\varepsilon}, x_{b})$. That is: further magnification by $\varepsilon$ of 
the snapshot (taken with magnification) $b$ is roughly the  same as the snapshot $b \varepsilon$. This 
is of course true in a neighbourhood of the base point $\displaystyle x_{b}$. 

The point (a) from the Definition \ref{dprofile}Ê has no other justification than Proposition \ref{propmetcone} in next 
subsection. 

We rewrite  definition \ref{defgh} with more details, in order to clearly understand what is a metric profile.  For any $b \in (0,a]$ and for any $\mu > 0$ there is $\varepsilon(\mu, b) \in (0,1)$ such that for any $\varepsilon \in 
(0,\varepsilon(\mu,b))$ there exists a relation $\displaystyle \rho = \rho_{\varepsilon, b} \subset \bar{B}_{d_{b}}(x_{b}, \varepsilon) \times \bar{B}_{d_{b \varepsilon}}(x_{b \varepsilon}, 1)$ such that
\begin{enumerate}
\itemsep-3pt
\item[1.] $\displaystyle dom \ \rho_{\varepsilon, b}$ is $\mu$-dense in $\displaystyle \bar{B}_{d_{b}}(x_{b}, \varepsilon)$, 
\item[2.] $\displaystyle im \ \rho_{\varepsilon, b}$ is $\mu$-dense in $\displaystyle \bar{B}_{d_{b \varepsilon}}(x_{b \varepsilon}, 1)$, 
\item[3.] $\displaystyle (x_{b}, x_{b\varepsilon}) \in \rho_{\varepsilon, b}$, 
\item[4.] for all $\displaystyle x,y \in \ dom \ \rho_{\varepsilon, b}$ we have
\begin{equation}
\label{disest}
\left|\frac{1}{\varepsilon} d_{b}(x,y) 
- d_{b \varepsilon}\left(\rho_{\varepsilon, b}(x), \rho_{\varepsilon, b}(y)\right)\right|
\, \leq \,  \mu
\end{equation}
\end{enumerate}

In the microscope interpretation, if $\displaystyle (x,u) \in \rho_{\varepsilon, b}$ means that $x$ and $u$ represent the same "real" point in the specimen. 

Therefore a metric profile gives two types of information:
\begin{itemize}
\itemsep-3pt
\item
a distance estimate like (\ref{disest}) from point 4,  
\item
an "approximate shape" estimate, like in the points 1--3, where we see that two sets, namely the balls   $\displaystyle \bar{B}_{d_{b}}(x_{b}, \varepsilon)$ and  $\displaystyle \bar{B}_{d_{b \varepsilon}}(x_{b \varepsilon}, 1)$, are approximately isometric. 
\end{itemize}

The simplest metric profile is one with $\displaystyle (\bar{B}(x_{b}, 1), d_{b}, x_{b}) = (X, d_{b}, x)$. 
In this case we see that $\displaystyle \rho_{\varepsilon, b}$ 
is approximately an $\varepsilon$ dilatation with base point $x$.  

This observation leads us to a particular class of (pointed) metric spaces, 
namely the  metric cones. 

\begin{definition}
\label{defmetcone}
A metric cone $(X,d, x)$ is a locally compact metric space $(X,d)$, with a marked point $x \in X$ such 
that for any $a,b \in (0,1]$ we have 
\[\displaystyle \mathbb{P}^{m}(a,x)  =  \mathbb{P}^{m}(b,x)\] 
\end{definition}

Metric cones have dilatations. By this we mean the following

\begin{definition}
Let $(X,d, x)$ be a metric cone. For any $\varepsilon \in (0,1]$  a dilatation is a function $\displaystyle \delta^{x}_{\varepsilon}: \bar{B}(x,1) \rightarrow \bar{B}(x,\varepsilon)$ such that 
\begin{itemize}
\itemsep-3pt
\item
$\displaystyle \delta^{x}_{\varepsilon}(x) = x$, 
\item
for any $u,v \in X$ we have 
\[d\left(\delta^{x}_{\varepsilon}(u), \delta^{x}_{\varepsilon}(v)\right) =\varepsilon \, d(u,v) \]
\end{itemize}
\end{definition}

The existence of dilatations for metric cones comes from the definition \ref{defmetcone}. 
Indeed, dilatations are just  isometries from $\displaystyle (\bar{B}(x,1), d, x)$ to $ (\bar{B}, \frac{1}{a}d, x)$. 

Metric cones are good candidates for being tangent spaces in the metric sense. 

\begin{definition}
\label{defmetspace}
A (locally compact) metric space $(M,d)$ admits a (metric) tangent space in $x \in M$ if the associated metric profile  
$\varepsilon \mapsto \mathbb{P}^{m}(\varepsilon, x)$ (as in definition \ref{dmprof})  admits a prolongation by continuity in 
$\varepsilon = 0$, i.e if the following limit exists: 
\begin{equation}
\label{limmetspace}
[T_{x}M,d^{x}, x]  =  \lim_{\varepsilon \rightarrow 0} \mathbb{P}^{m}(\varepsilon, x)
\end{equation}
\end{definition}

The connection between metric cones, tangent spaces and metric profiles in the abstract sense is made by the following proposition. 

\begin{proposition}
\label{propmetcone}
The associated metric profile $\varepsilon \mapsto \mathbb{P}^{m}(\varepsilon, x)$ of a metric space $(M,d)$ for  a fixed $x \in M$ is a metric profile in the sense of the definition \ref{dprofile} if and only if the space 
$(M,d)$ admits a tangent space in $x$. 
In such a case the tangent space is a metric cone. 
\end{proposition}

\begin{proof} 
A tangent space  $[V,d_{v}, v]$ exists if and only if we have the limit from the relation (\ref{limmetspace}). 
In this case there exists a prolongation by continuity to $\varepsilon = 0$  of the  metric profile 
$\mathbb{P}^{m}(\cdot , x)$. 
The prolongation is a metric profile in the sense of definition \ref{dprofile}. 
Indeed, we have still to check the property (b). But this is trivial, because for any $\varepsilon, b >0$,  sufficiently small, we have 
\[
\mathbb{P}^{m}(\varepsilon b, x) =  \mathbb{P}^{m}_{d_{b}}(\varepsilon,x)
\]
where  $d_{b} = (1/b) d$ and 
$\mathbb{P}^{m}_{d_{b}}(\varepsilon,x) = [\bar{B}(x,1),\frac{1}{\varepsilon}d_{b}, x]$.

Finally, let us prove that the tangent space is a metric cone. For any $a \in (0,1]$ we have
\[
\left[\bar{B}(x,1), \frac{1}{a}d^{x}, x\right] = \lim_{\varepsilon \rightarrow 0} \mathbb{P}^{m}(a \varepsilon, x) \]
Therefore 
\[
\left[\bar{B}(x,1), \frac{1}{a}d^{x}, x\right]  = [T_{x}M,d^{x}, x] 
 \tag*{\qed}
\] 
\renewcommand{\qed}{}
\end{proof}

\subsection{Groups with dilatations. Virtual tangent space}

In section 6 we shall see  that metric tangent spaces sometimes have a group structure which is compatible with 
dilatations. This structure, of a group with dilatations, is interesting by itself. The notion has been introduced in \cite{buliga2}; we describe  it further. 

We start with the following setting: $G$ is a topological group endowed with an uniformity such that the operation is uniformly continuous.  The description that follows  is slightly non canonical, but is nevertheless 
motivated by the case of a Lie group endowed with a Carnot-Caratheodory  distance induced by a left invariant distribution. 

We introduce first the double of $G$, as the group $G^{(2)} = G \times G$ with operation
\[(x,u) (y,v) = (xy, y^{-1}uyv) \]
The operation on the group $G$, seen as the function
\[op: G^{(2)} \rightarrow G \ , \ \ op(x,y) = xy \]
is a group morphism. Also the inclusions:
\begin{gather*}
i': G \rightarrow G^{(2)} , \quad i'(x) = (x,e)\\
i": G \rightarrow G^{(2)} , \quad i"(x) = (x,x^{-1})
\end{gather*}
are group morphisms.

\begin{definition}
\label{dunifg}
\begin{enumerate}
\itemsep-3pt
\item[1.]
$G$ is a uniform group if we have two uniformity structures, on $G$ and
$G^{2}$,  such that $op$, $i'$, $i"$ are uniformly continuous.
\item[2.] 
A local action of a uniform group $G$ on a uniform  pointed space $(X, x_{0})$ is a function
$\phi \in W \in \mathcal{V}(e)  \mapsto \hat{\phi}: U_{\phi} \in \mathcal{V}(x_{0}) \rightarrow
V_{\phi}  \in \mathcal{V}(x_{0})$ such that
\begin{enumerate}
\itemsep-3pt
\item[(a)] 
the map $(\phi, x) \mapsto \hat{\phi}(x)$ is uniformly continuous from $G \times X$ (with product uniformity)
to  $X$,
\item[(b)] 
for any $\phi, \psi \in G$ there is $D \in \mathcal{V}(x_{0})$
such that for any $x \in D$ $\hat{\phi \psi^{-1}}(x)$ and $\hat{\phi}(\hat{\psi}^{-1}(x))$ make sense and   $\hat{\phi \psi^{-1}}(x) = \hat{\phi}(\hat{\psi}^{-1}(x))$.
\end{enumerate}
\item[3.] 
Finally, a local group is an uniform space $G$ with an operation defined in a neighbourhood of $(e,e) \subset G \times G$ which satisfies the uniform group axioms locally.
\end{enumerate}
\end{definition}
Note that a local group acts locally at left (and also by conjugation) on itself.

This definition deserves an explanation. An uniform group, according to the Definition \ref{dunifg}, is a group $G$ 
such that left translations are uniformly continuous functions and the left action of $G$ on itself is 
uniformly continuous too. In order to precisely formulate this we need two uniformities: one on $G$ and another 
on $G \times G$. 

These uniformities should be compatible, which is achieved by saying that $i'$, $i"$ are uniformly continuous. 
The uniformity of the group operation is achieved by saying that the $op$ morphism is uniformly continuous. 

\begin{definition}
\label{defgwd}
A group with dilatations $(G,\delta)$ is a local uniform group $G$  with  a local action of $\Gamma$ (denoted by $\delta$), on $G$ such that
\begin{enumerate}
\itemsep-3pt
\item[H0.] the limit  $\displaystyle \lim_{\varepsilon \rightarrow 0} \delta_{\varepsilon} x = e$ exists and is uniform with respect to $x$ in a compact neighbourhood of the identity $e$.
\item[H1.] the limit
\[\beta(x,y) = \lim_{\varepsilon \rightarrow 0} \delta_{\varepsilon}^{-1}
\left((\delta_{\varepsilon}x) (\delta_{\varepsilon}y ) \right)\]
is well defined in a compact neighbourhood of $e$ and the limit is uniform.
\item[H2.] the following relation holds:
\[ \lim_{\varepsilon \rightarrow 0} \delta_{\varepsilon}^{-1}
\left( ( \delta_{\varepsilon}x)^{-1}\right) = x^{-1}\]
where the limit from the left hand side exists in a neighbourhood of $e$ and is uniform with respect to $x$.
\end{enumerate}
\end{definition}

These axioms are the prototype of a dilatation structure. 

 The "infinitesimal version" of an uniform group is a conical 
local uniform group. 

\begin{definition}
A conical group $N$ is a local group with a local action of
$(0,+\infty)$ by morphisms $\delta_{\varepsilon}$ such that
$\displaystyle \lim_{\varepsilon \rightarrow 0} \delta_{\varepsilon} x = e$ for any
$x$ in a neighbourhood of the neutral element $e$.
\end{definition}

Here comes a proposition which explains why a conical group is the infinitesimal version of a group with 
dilatations. 

\begin{proposition}
\label{here3.4}
Under the hypotheses H0, H1, H2 $(G,\beta, \delta)$ is a conical group, with operation 
$\beta$ and dilatations $\delta$.
\end{proposition}

\begin{proof}
All the uniformity assumptions allow us  to change at will the order of taking
limits. We shall not insist on this further and we shall concentrate on the
algebraic aspects.

We have to prove the associativity, existence of neutral element, existence of 
inverse and the property of being conical. 

For the associativity $\beta(x,\beta(y,z))  =  \beta(\beta(x,y),z)$ we calculate
\[\beta(x,\beta(y,z))  =  \lim_{\varepsilon \rightarrow 0 , \eta \rightarrow 0} \delta_{\varepsilon}^{-1} 
\left\{ (\delta_{\varepsilon}x) \delta_{\varepsilon/\eta}\left( (\delta_{\eta}y) (\delta_{\eta} z) \right) \right\} \]
We take $\varepsilon = \eta$ and get
\[ \beta(x,\beta(y,z))  =  \lim_{\varepsilon \rightarrow 0}\left\{
(\delta_{\varepsilon}x) (\delta_{\varepsilon} y) (\delta_{\varepsilon} z) \right\} \]
In the same way
\[\beta(\beta(x,y),z)  =  \lim_{\varepsilon \rightarrow 0 , \eta \rightarrow 0} \delta_{\varepsilon}^{-1} 
\left\{ (\delta_{\varepsilon/\eta}x)\left( (\delta_{\eta}x) (\delta_{\eta} y) \right) (\delta_{\varepsilon} z) 
\right\} \]
and again taking $\varepsilon = \eta$ we obtain
\[\beta(\beta(x,y),z)  =   \lim_{\varepsilon \rightarrow 0}\left\{
(\delta_{\varepsilon}x) (\delta_{\varepsilon} y) (\delta_{\varepsilon} z) \right\} = \beta(x,\beta(y,z)) \]
The neutral element is $e$, from H0 (first part) it follows that $\beta(x,e)  = \beta(e,x)  =  x$. 
The inverse of $x$ is $x^{-1}$, by a similar argument:
\[\beta(x, x^{-1})   =  \lim_{\varepsilon \rightarrow 0 , \eta \rightarrow 0} 
\delta_{\varepsilon}^{-1} \left\{ (\delta_{\varepsilon}x)
\left( \delta_{\varepsilon/\eta}(\delta_{\eta}x)^{-1}\right) \right\} \]
and taking $\varepsilon = \eta$ we obtain
\[\beta(x, x^{-1})   =  \lim_{\varepsilon \rightarrow 0}
\delta_{\varepsilon}^{-1} \left( (\delta_{\varepsilon}x) (\delta_{\varepsilon}x)^{-1}\right) 
 =  \lim_{\varepsilon \rightarrow 0} \delta_{\varepsilon}^{-1}(e)  = \ e \]
Finally, $\beta$ has the property
\[\beta(\delta_{\eta} x, \delta_{\eta}y)  =  \delta_{\eta} \beta(x,y) \]
which comes from the definition of $\beta$ and commutativity of multiplication
in $(0,+\infty)$. This proves that $(G,\beta, \delta)$ is conical. 
\end{proof}

In a sense $(G,\beta,\delta)$ is the tangent space of the group with dilatations $(G,\delta)$ at $e$. We can 
act with  the conical group  $(G,\beta,\delta)$ on $(G,\delta)$. Indeed, let us denote by $[f,g] = f \circ g \circ f^{-1} \circ g^{-1}$ the commutator of two transformations. For the group $G$  we shall denote by
$L_{x}^{G} y = xy$ the left translation and by $L^{N}_{x}y = \beta(x,y)$. The preceding proposition tells us that $(G,\beta, \delta)$ acts locally by left
translations on $G$. We shall call the left translations with respect to the group operation $\beta$ "infinitesimal". 
These infinitesimal translations admit an interesting commutator representation
\begin{equation}
\lim_{\lambda \rightarrow 0}\left[L_{(\delta_{\lambda}x)^{-1}}^{G}, \delta_{\lambda}^{-1}\right] = L^{N}_{x} 
\label{firstdef}
\end{equation}

\begin{definition}
The group $VT_{e}G$ formed by all transformations $L_{x}^{N}$ is called the virtual tangent space at $e$ to $G$.
\end{definition}

As local groups,  $VT_{e}G$ and $(G,\beta, \delta)$ are isomorphic. We can easily define dilatations on  $VT_{e}G$, by conjugation with dilatations $\delta_{\varepsilon}$. Really, we see that 
\[
L^{N}_{\delta_{\varepsilon}x}(y) 
= \beta(\delta_{\varepsilon} x, y)
= \delta_{\varepsilon} L^{N}_{x} \left( \delta_{\varepsilon} \right)^{-1}
\]

The virtual tangent space $VT_{x}G$ at $x \in G$ to $G$ is obtained by translating the group operation and the 
dilatations from $e$ to $x$. This means: define a new operation on $G$ by
\[y \stackrel{x}{\cdot} z = y x^{-1}z \]
The group $G$ with this operation is isomorphic to $G$ with old operation and
the left translation $L^{G}_{x}y = xy$ is the isomorphism. The neutral element is $x$.
Introduce also the dilatations based at $x$ by
\[\delta_{\varepsilon}^{x} y = x \delta_{\varepsilon}(x^{-1}y) \]
Then $G^{x} = (G,\stackrel{x}{\cdot})$ with the group of dilatations $\delta_{\varepsilon}^{x}$ satisfy the Axioms 
H0, H1, H2. Define then the virtual tangent
space $VT_{x}G$ to be: $VT_{x}G = VT_{x} G^{x}$. 

\section{Binary decorated trees and dilatations}

We want to explore what happens when we make compositions of dilatations 
(which depends also on $\varepsilon> 0$ ). The $\varepsilon$ variable apart, any dilatation 
$\displaystyle \delta_{\varepsilon}^{x}(y)$ is a function of two arguments: 
$x$ and $y$, invertible with respect to the second argument. 
The functions we can obtain when composing dilatations are 
difficult to write, that is why we shall use a tree notation. 

\subsection{The formalism}

Let $X$ be a non empty set and $\mathcal{T}(X)$ be a  class
 of binary planar trees with  leaves in $X$ and all  nodes  
 decorated with two colors $\displaystyle \left\{ \circ , 
\bullet \right\}$. The empty tree, that is the tree with no nodes or leaves, 
belongs to $\mathcal{T}(X)$.  
For any $x\in X$ we accept that there is a tree in $\mathcal{T}(X)$ 
 with no nodes and with $x$ as the only leaf. 
That is $X \subset \mathcal{T}(X)$. 

For any color $\mathbf{a} \in \displaystyle 
\left\{ \circ , \bullet \right\}$, let $\bar{\mathbf{a}}$ be the 
opposite color. The colors $\circ$ and $\bullet$ are codes for the symbols 
$\varepsilon$ and $\displaystyle \varepsilon^{-1}$.

The relation "$\approx$" is an 
equivalence relation on $\mathcal{T}(X)$, taken as a primitive notion 
for the axioms which will follow. 

The equivalence class of of a tree $\mathcal{P} \in \mathcal{T}(X)$ is denoted 
by $\displaystyle\begin{parsetree} (...$\mathcal{P}$.)\end{parsetree}$. 
In various diagrams that will follow we shall use the notation 
$\displaystyle  \begin{parsetree}(.$\Gamma$. .$\mathcal{P}$.)\end{parsetree}$ 
for saying that $\Gamma$ is the equivalence class of $\mathcal{P}$. 
For any 
$\displaystyle \mathcal{P}, \mathcal{R} \in
\mathcal{T}(X)$,  "$\mathcal{P}\approx \mathcal{R}$"   or
"$\begin{parsetree} (.. .$\mathcal{P}$.)\end{parsetree} 
=  \begin{parsetree} (.. .$\mathcal{R}$.)\end{parsetree}$"  
means  the same thing. 

\bigskip
{\bf Axiom T0.} 
For any $x, u, v \in X$ the trees 
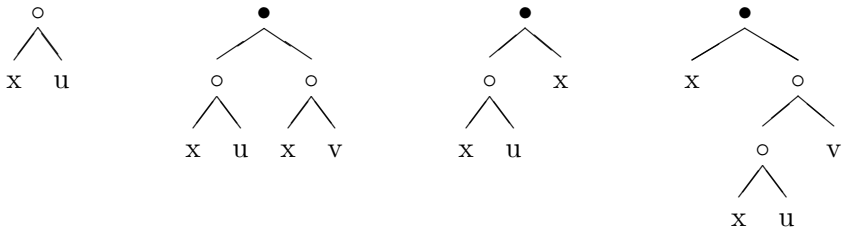
\[ 
\begin{parsetree}
    ( .$\circ$.  .x. .u.)
\end{parsetree} 
\qquad\qquad 
\begin{parsetree}
    ( .$\bullet$.  
                                      (.$\circ$. .x. .u.)  (.$\circ$. .x. .v.)
    )
\end{parsetree} 
\qquad\qquad  
\begin{parsetree}
    ( .$\bullet$.  
                                      (.$\circ$. .x. .u.)  .x.
    )
\end{parsetree} 
\qquad\qquad  
\begin{parsetree}
    ( .$\bullet$.  
                                      .x.  
                                             (.$\circ$. (.$\circ$. .x. .u.)  .v.)
    )\end{parsetree} 
\]
belong to $\mathcal{T}(X)$.

The equivalence class of 
$\begin{parsetree} ( .$\circ$.  .x. .u.)\end{parsetree}$ 
is denoted by $\displaystyle \delta^{x}_{\varepsilon} u$, that is we have
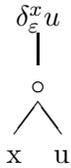
\[
\begin{parsetree}
 (.$\delta^{x}_{\varepsilon} u$.
    ( .$\circ$.  .x. .u.))
\end{parsetree}
\] 

{\bf Axiom T1.} 
Consider any  trees $\mathcal{P}, \mathcal{R}, \mathcal{S}, 
\mathcal{Q}, \mathcal{Z}  \in \mathcal{T}(X)$, any $x \in X$,  
and any colors $\mathbf{a}, \mathbf{b}$  such that the trees from the right 
hand sides of relations below belong to $\mathcal{T}(X)$. Then  the trees 
 from the left hand sides of relations below belong to $\mathcal{T}(X)$ 
and we have
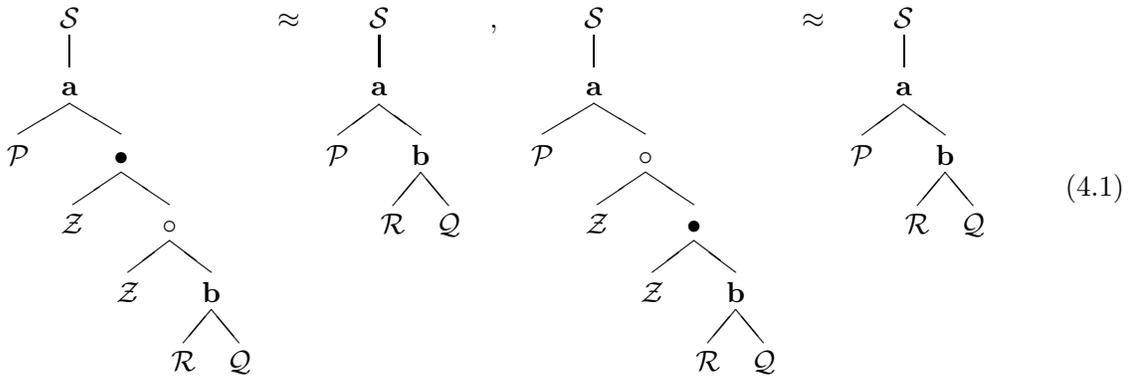
\begin{equation}
\begin{array}{ccccccc}
\begin{parsetree}
(.$\mathcal{S}$. (.$\mathbf{a}$. .$\mathcal{P}$.
    ( .$\bullet$.  
                     .$\mathcal{Z}$.  (.$\circ$. .$\mathcal{Z}$. 
(.$\mathbf{b}$. .$\mathcal{R}$. 
.$\mathcal{Q}$.  ))
    )))
\end{parsetree} &  \approx  & \begin{parsetree} 
(.$\mathcal{S}$.  (.$\mathbf{a}$. .$\mathcal{P}$.
     (.$\mathbf{b}$. .$\mathcal{R}$. 
.$\mathcal{Q}$.  )
    ))
 \end{parsetree} 
\quad  , \quad
\begin{parsetree}
(.$\mathcal{S}$. (.$\mathbf{a}$. .$\mathcal{P}$.
    ( .$\circ$.  
                     .$\mathcal{Z}$.  (.$\bullet$. .$\mathcal{Z}$. (.$\mathbf{b}$. .$\mathcal{R}$. 
.$\mathcal{Q}$.  ))
    )))
\end{parsetree} &  \approx  & \begin{parsetree} 
(.$\mathcal{S}$.  (.$\mathbf{a}$. .$\mathcal{P}$.
     (.$\mathbf{b}$. .$\mathcal{R}$. 
.$\mathcal{Q}$.  )
    ))
 \end{parsetree}
\end{array}
\label{ax1}
\end{equation}
Here, in all  diagrams,  the symbol $\displaystyle \begin{parsetree} 
(.$\mathcal{S}$.  .$\mathbf{a}$.)  \end{parsetree}$  means that the node 
colored with $\mathbf{a}$ is grafted at an arbitrary leaf of the 
tree $\mathcal{S}$.

The second axiom expresses the fact that the dilatation (of any coefficient 
$\varepsilon$) $\displaystyle \delta^{x}_{\varepsilon}$ has $x$ as fixed point, 
that is $\displaystyle  \delta^{x}_{\varepsilon} x = x$. 

{\bf Axiom T2.} 
For any $x\in X$ the tree $\begin{parsetree}( .$\bullet$.  .x. .x.)\end{parsetree}$ 
belongs to $\mathcal{T}(X)$. Moreover, consider any  tree 
$\mathcal{P} \in \mathcal{T}(X)$ and  any $x \in X$.  Then  the trees 
from the left hand sides of relations below belong to $\mathcal{T}(X)$ and we have
\begin{equation}
\begin{array}{ccccccc}
 \begin{parsetree}
 (.$\mathcal{P}$.
    ( .$\circ$.  .x. .x.)) 
\end{parsetree} & \approx & \mathcal{P},
\qquad 
\begin{parsetree}
 (.$\mathcal{P}$.
    ( .$\bullet$.  .x. .x.)) 
\end{parsetree} & \approx & \mathcal{P} 
\end{array}
\label{ax2}
\end{equation}
that is the equivalence class of $x$ is the same as the equivalence class of 
$\begin{parsetree}
    ( .$\circ$.  .x. .x.)
\end{parsetree}$ and the equivalence class of $\begin{parsetree}
    ( .$\bullet$.  .x. .x.)
\end{parsetree}$. As in Axiom T1, the symbol $\displaystyle \begin{parsetree} 
(.$\mathcal{S}$.  .$\mathcal{P}$.)  \end{parsetree}$  means that the root of the tree $\mathcal{P}$  
 is grafted at an arbitrary leaf of the 
tree $\mathcal{S}$.

\begin{definition}
\label{defsumdif}
We define the difference, sum and inverse trees as follows:
\begin{enumerate}
\itemsep-3pt
\item[(a)] 
the difference tree $\displaystyle \Delta^{x}_{\varepsilon} = 
\Delta^{x}_{\varepsilon}(u,v)$ is given by the relation 
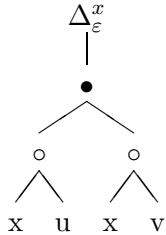
\[
\begin{parsetree}
 (.$\Delta^{x}_{\varepsilon}$.
    ( .$\bullet$.  
                                      (.$\circ$. .x. .u.)  (.$\circ$. .x. .v.)
    ))
\end{parsetree}
\]
\item[(b)] 
the sum tree $\displaystyle \Sigma^{x}_{\varepsilon} = \Sigma^{x}_{\varepsilon}(u,v)$ is
given by the relation 
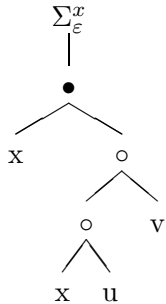
\[\begin{parsetree}
 (.$\Sigma^{x}_{\varepsilon}$.
    ( .$\bullet$.  
                                      .x.  
                                             (.$\circ$. (.$\circ$. .x. .u.)  .v.)
    ))\end{parsetree}\]
\item[(c)] 
the inverse tree $\displaystyle inv^{x}_{\varepsilon} = inv^{x}_{\varepsilon}(u)$ 
is given by the relation 
    \[\begin{parsetree}
 (.$inv^{x}_{\varepsilon}$.
    ( .$\bullet$.  
                                      (.$\circ$. .x. .u.)  .x.
    ))
\end{parsetree}
\]
\end{enumerate}
\end{definition}

The next axiom states that T0, T1, T2 are sufficient for determining the class 
$\mathcal{T}(X)$ and the equivalence relation $\approx$. 

{\bf Axiom T3.} 
The class $\mathcal{T}(X)$ is the smallest class of trees 
obtained by grafting of trees listed in Axiom T0, and satisfying Axioms T1, T2.
Moreover, two trees from $\mathcal{T}(X)$ are equivalent if and only if they 
can be proved equivalent after a finite string of applications of Axioms T1, T2.

\subsection{First consequences}

We shall use  the axioms in order to obtain results that we shall use later, for 
dilatation structures. 

\begin{proposition}
\label{psumdif}
For any $x, u, y$ and $v$ we have 
\begin{enumerate}
\itemsep-3pt
\item[(a)] $ \displaystyle \Delta^{x}_{\varepsilon}(u, \Sigma^{x}_{\varepsilon}(u,y)) = y$, 
\item[(b)] $\displaystyle \Sigma^{x}_{\varepsilon}(u,\Delta^{x}_{\varepsilon}(u,v)) = v$.
\end{enumerate}
\end{proposition}

\begin{proof}
We prove (a) by computations using the definition \ref{defsumdif} of the sum and 
difference trees, and Axiom T1 several times.   
\[\begin{array}{ccccc}
\begin{parsetree}
( .. ( .$\bullet$. ( .$\circ$. .x. .u.) ( .$\circ$. .x. ( .$\bullet$. .x. ( .$\circ$. ( .$\circ$. .x. .u.  ) .y.  ) ))     ))
\end{parsetree} & = & \begin{parsetree}
( .. ( .$\bullet$. ( .$\circ$. .x. .u.) ( .$\circ$. ( .$\circ$. .x. .u.   )   .y.        )          )) 
\end{parsetree} & = & y  
\end{array}\]
For (b) we proceed in the same way:
\[
\begin{array}{ccccc}
\begin{parsetree}
( .. ( .$\bullet$. .x. ( .$\circ$. ( .$\circ$. .x. .u.) ( .$\bullet$. ( .$\circ$. .x. .u.)  ( .$\circ$. .x. .v.)    )  )  )) 
\end{parsetree} & = & \begin{parsetree}
( .. ( .$\bullet$. .x. ( .$\circ$. .x. .v.) ))
\end{parsetree} & = & v      
\end{array}
 \tag*{\qed}
\]
\renewcommand{\qed}{}
\end{proof}

\begin{proposition}
\label{p33}
We have the relations
\begin{equation}
\label{eqinv}
 \Delta^{x}_{\varepsilon}(u,v)  = \Sigma^{ \quad \begin{parsetree}
    ( .$\circ$.  .x. .u.)
\end{parsetree} \quad }_{\varepsilon}(inv^{x}_{\varepsilon}(u),v) 
\end{equation}
\begin{equation}
\label{difinv}
 inv^{x}_{\varepsilon}(u)  = \Delta^{x}_{\varepsilon}(u,x) 
\end{equation}
\begin{equation}
inv^{ \quad \begin{parsetree}
    ( .$\circ$.  .x. .u.)
\end{parsetree} \quad }_{\varepsilon}\left( (inv^{x}_{\varepsilon}(u) \right) = u
\label{invinv}
\end{equation}
\end{proposition} 

\begin{proof}
Graphically, the relation (\ref{eqinv}) is
\[
\begin{array}{ccc}
\begin{parsetree}
 (..
    ( .$\bullet$.  
                                      (.$\circ$. .x. .u.)  (.$\circ$. .x. .v.)
    ))
\end{parsetree} & = & \begin{parsetree}
 (..
    ( .$\bullet$.  
                                      ( .$\circ$.  .x. .u.)  
                                             (.$\circ$. (.$\circ$. ( .$\circ$.  .x. .u.) ( .$\bullet$.  
                                      (.$\circ$. .x. .u.)  .x.
    ))  .v.)
    ))\end{parsetree}
 \end{array} 
\]
This is true by Axiom T1. 

The relation (\ref{difinv}) is 
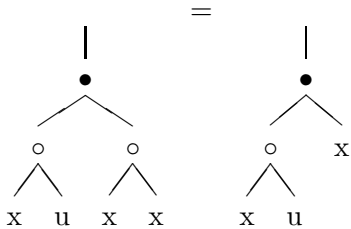
\[
\begin{array}{ccc}
\begin{parsetree}
 (..
    ( .$\bullet$.  
                                      (.$\circ$. .x. .u.)  (.$\circ$. .x. .x.)
    ))
\end{parsetree} & = & \begin{parsetree}
 (..
    ( .$\bullet$.  
                                      (.$\circ$. .x. .u.)  .x.
    ))
\end{parsetree}
 \end{array}
 \]
This is true by Axiom T2. 

We prove the relation  (\ref{invinv}) by a string of equalities, starting from the left hand side to the right: 
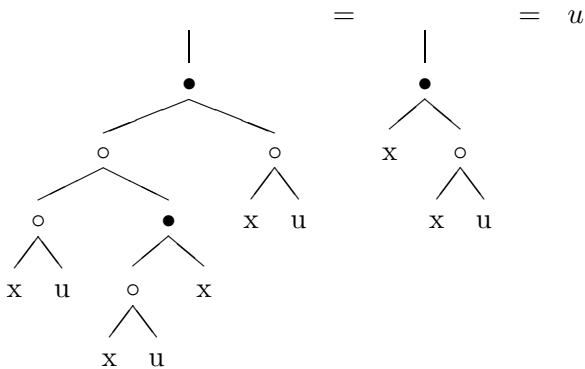
\[\begin{array}{ccccc}
\begin{parsetree}
 (..
    ( .$\bullet$.  
                                      (.$\circ$. ( .$\circ$.  .x. .u.) ( .$\bullet$.  
                                      (.$\circ$. .x. .u.)  .x.
    ))  ( .$\circ$.  .x. .u.)
    ))
\end{parsetree} & = & \begin{parsetree}
 (..
    ( .$\bullet$.  
                                      .x.   (.$\circ$. .x. .u.)
    ))
    \end{parsetree} & = & u
\end{array} \]
Here we have used the Axiom T1 several times. 
\end{proof}

The relation (\ref{invinv}) in last proposition shows that the "inverse function" $\displaystyle inv^{x}_{\varepsilon}$ is not involutive, but shifted involutive. 

The next proposition shows that the function $\displaystyle  \Sigma^{x}_{\varepsilon}(\cdot,\cdot)$  satisfies 
a shifted associativity property. 

\begin{proposition}
\label{p34}
We have the following relations:
\begin{equation}
\label{eqasoc1}
 \Delta^{x}_{\varepsilon}(u,\Sigma^{x}_{\varepsilon}(\Sigma^{x}_{\varepsilon}(u,v),w))  
= \Sigma^{\quad\begin{parsetree}
    ( .$\circ$.  .x. .u.)
\end{parsetree} \quad }_{\varepsilon}(v,w) 
\end{equation}

\begin{equation}
\label{eqasoc2}
 \Sigma^{x}_{\varepsilon}\left(u,\Sigma^{ \quad \begin{parsetree}
 ( .$\circ$.  .x. .u.)
\end{parsetree} \quad }_{\varepsilon}(v,w)\right)  = \Sigma^{x}_{\varepsilon}(\Sigma^{x}_{\varepsilon}(u,v),w)
\end{equation}
\end{proposition} 

\begin{proof}
Graphically, the relation (\ref{eqasoc1}) is
\[\begin{array}{ccc}
\begin{parsetree}
 (..
    ( .$\bullet$.  
                                      (.$\circ$. .x. .u.)  (.$\circ$. .x. ( .$\bullet$.  
                                      .x.  
                                             (.$\circ$. (.$\circ$. .x. ( .$\bullet$.  
                                      .x.  
                                             (.$\circ$. (.$\circ$. .x. .u.)  .v.)
    ))  .w.)
    )  )
    )) 
\end{parsetree} 
& = & 
\begin{parsetree}
(.. ( .$\bullet$.  
                                       ( .$\circ$.  .x. .u.)
                                             (.$\circ$. (.$\circ$.  
    ( .$\circ$.  .x. .u.) .v.)  .w.)
    ))
 \end{parsetree}
\end{array}\]

This is true by Axiom T1. 

The relation (\ref{eqasoc2}) is is equivalent to (\ref{eqasoc1}), by Proposition \ref{psumdif}.  We can also give a direct proof by graphically representing the relation
\[\begin{array}{ccc}
\begin{parsetree}
 (..
    ( .$\bullet$.  
                                      .x.  
                                             (.$\circ$. (.$\circ$. .x. .u.)  ( .$\bullet$.  
                                       ( .$\circ$.  .x. .u.)
                                             (.$\circ$. (.$\circ$.  
    ( .$\circ$.  .x. .u.) .v.)  .w.)
    ))
    ))\end{parsetree} & = & \begin{parsetree}
 (..
    ( .$\bullet$.  
                                      .x.  
                                             (.$\circ$. (.$\circ$. .x. ( .$\bullet$.  
                                      .x.  
                                             (.$\circ$. (.$\circ$. .x. .u.)  .v.)
    ))  .w.)
    ))\end{parsetree}
    \end{array} \]
This is true by the Axiom T1. 
\end{proof}

\section{Dilatation structures}

The space $(X,d)$ is a complete, locally compact metric space. This means that as a metric space 
$(X,d)$ is complete and  that small balls are compact. 

\subsection{Axioms of dilatation structures}

The  axioms of  a dilatation structure $(X,d,\delta)$ are listed further. 
The first axiom is merely a preparation for the next axioms. That is why we 
counted it as Axiom 0.

\begin{enumerate}
\itemsep-3pt
\item[{\bf A0.}] 
Depending on the parameter $\varepsilon \in (0,+\infty)$, dilatations are objects 
having the following description. 

For any $\displaystyle \varepsilon \in (0,1]$ 
the dilatations are functions \[ \delta_{\varepsilon}^{x}: U(x) \rightarrow 
V_{\varepsilon}(x) \]
All  such dilatations  are homeomorphisms (invertible, continuous, with continuous inverse). 

We suppose  that there is  $1<A$ such that for any $x \in X$ we have 
\[\bar{B}_{d}(x,A) \subset U(x)  \]
 We suppose that for all  $\varepsilon \in 
(0,1)$, we have 
\[ B_{d}(x,\varepsilon) \subset \delta_{\varepsilon}^{x} B_{d}(x,A) \subset V_{\varepsilon}(x) \subset U(x) \]

 For $\varepsilon \in (1,+\infty)$ the associated dilatation  
\[\delta^{x}_{\varepsilon} : W_{\varepsilon}(x) \rightarrow B_{d}(x,B) \ , \]
is an injective, continuous, with continuous inverse  on the image. 
We shall suppose that $\displaystyle  W_{\varepsilon}(x)$ is open, 
\[V_{\varepsilon^{-1}}(x) \subset W_{\varepsilon}(x)\]
and that for all $\displaystyle \varepsilon \in [0,1]$ and 
$\displaystyle u \in U(x)$ we have
\[\delta_{\varepsilon^{-1}}^{x} \ \delta^{x}_{\varepsilon} u = u \]
\end{enumerate}

We remark that we have the following string of inclusions, for any 
$\varepsilon \in (0,1]$  and any $x \in X$:
\[ B_{d}(x,\varepsilon) \subset \delta^{x}_{\varepsilon}  B_{d}(x, A) \subset V_{\varepsilon}(x) \subset 
W_{\varepsilon^{-1}}(x) \subset \delta_{\varepsilon}^{x}  B_{d}(x, B) \]

A further technical condition on the sets  $\displaystyle V_{\varepsilon}(x)$ and $\displaystyle W_{\varepsilon}(x)$  will be given just before the Axiom A4. (This condition will be counted as part of 
Axiom A0.)

\begin{enumerate}
\itemsep-3pt
\item[{\bf A1.}]  
We  have $\displaystyle  \delta^{x}_{\varepsilon} x = x $ for any point $x$. 
We also have $\displaystyle \delta^{x}_{1} = id$ for any $x \in X$.

Let us define the topological space
\[ dom \, \delta = \left\{ (\varepsilon, x, y) \in (0,\infty) \times X 
\times X \mbox{: } \mbox{if } \varepsilon \in (0,1] \mbox{ then } y \in U(x) \mbox{, else } y \in W_{\varepsilon}(x) \right\} \] 
with the topology inherited from the product topology on 
$\Gamma \times X \times X$. Consider also 
$\displaystyle Cl(dom \, \delta)$, the closure of $dom \, \delta$ in 
$\displaystyle [0,\infty) \times X \times X$ with product topology. 
The function 
\[\delta : dom \, \delta \rightarrow  X\] defined by 
$\displaystyle \delta (\varepsilon,  x, y)  = \delta^{x}_{\varepsilon} y$ is continuous. Moreover, it can be continuously extended to $\displaystyle Cl(dom \, \delta)$ and we have 
\[\lim_{\varepsilon\rightarrow 0} \delta_{\varepsilon}^{x} y \, = \, x \]
\item[{\bf A2.}] 
For any  $x, \in K$, $\displaystyle \varepsilon, \mu \in \Gamma_{1}$ and $\displaystyle u \in 
\bar{B}_{d}(x,A)$   we have
\[ \delta_{\varepsilon}^{x} \delta_{\mu}^{x} u  = \delta_{\varepsilon \mu}^{x} u  \]
\item[{\bf A3.}]
 For any $x$ there is a  function 
$\displaystyle (u,v) \mapsto d^{x}(u,v)$, defined for any $u,v$ in the closed ball (in distance d) 
$\displaystyle\bar{B}(x,A)$, such that 
\[\lim_{\varepsilon \rightarrow 0} \quad \sup  \left\{\left|\frac{1}{\varepsilon} d(\delta^{x}_{\varepsilon} u, \delta^{x}_{\varepsilon} v) - d^{x}(u,v) \right|\mbox{:  } u,v \in \bar{B}_{d}(x,A)\right\} =0\]
uniformly with respect to $x$ in compact set. 
\end{enumerate}

\begin{rk}
\label{imprk}
The "distance" $d^{x}$ can be degenerated. That means: there might be 
$\displaystyle v,w \in \bar{B}_{d}(x,A)$ such that $\displaystyle d^{x}(v,w) = 0$ but $v \not = w$. We shall 
use further the name "distance" for $d^{x}$, essentially by commodity, but  keep 
in mind the possible degeneracy of $d^{x}$. 
\end{rk}

For  the following axiom to make sense we impose a technical condition on the co-domains $\displaystyle V_{\varepsilon}(x)$: for any compact set $K \subset X$ there are $R=R(K) > 0$ and 
$\displaystyle \varepsilon_{0}= \varepsilon(K) \in (0,1)$  such that  
for all $\displaystyle u,v \in \bar{B}_{d}(x,R)$ and all $\displaystyle \varepsilon \in \Gamma$, $\displaystyle  \nu(\varepsilon) \in (0,\varepsilon_{0})$,  we have 
\[\delta_{\varepsilon}^{x} v \in W_{\varepsilon^{-1}}( \delta^{x}_{\varepsilon}u) \]

With this assumption the following notation makes sense:
\begin{equation}
\label{delnot}
\Delta^{x}_{\varepsilon}(u,v) = 
\delta_{\varepsilon^{-1}}^{\delta^{x}_{\varepsilon} u} 
\delta^{x}_{\varepsilon} v
\end{equation}
The next axiom can now be stated: 
\begin{enumerate}
\itemsep-3pt
\item[{\bf A4.}] 
We have the limit 
\[\lim_{\varepsilon \rightarrow 0}  \Delta^{x}_{\varepsilon}(u,v) =  \Delta^{x}(u, v)  \]
uniformly with respect to $x, u, v$ in compact set. 
\end{enumerate}

Note that with the tree notation we may identify (\ref{delnot}) with the difference tree  from 
Definition \ref{defsumdif} (a).  

\begin{definition}
A triple $(X,d,\delta)$ which satisfies A0, A1, A2, A3, but $\displaystyle d^{x}$ is degenerate for some 
$x\in X$, is called degenerate dilatation structure. 

If the triple $(X,d,\delta)$ satisfies A0, A1, A2, A3 and 
 $\displaystyle d^{x}$ is non-degenerate for any $x\in X$, then we call it  a weak dilatation structure. 

 If a weak dilatation structure satisfies A4 then we call it 
dilatation structure. 
 \end{definition}

Note that it could be assumed, without great modification of the axioms, that
\begin{enumerate}
\itemsep-3pt
\item[(a)] 
we may replace $(0,\infty)$ by  $\Gamma$,  
a topological separated commutative group  endowed with a continuous group 
morphism $\displaystyle \nu : \Gamma \rightarrow (0,+\infty)$ with 
$\displaystyle \inf \nu(\Gamma)  =  0$. Here $(0,+\infty)$ is 
taken as a group with multiplication. The neutral element of $\Gamma$ is 
denoted by $1$. We use the multiplicative notation for the operation in $\Gamma$. 

The morphism $\nu$ defines an invariant topological filter on $\Gamma$ (equivalently, an end). Really, 
this is the filter generated by the open sets $\displaystyle \nu^{-1}(0,a)$, $a>0$. From now on 
we shall name this topological filter (end) by "0" and we shall write $\varepsilon \in \Gamma \rightarrow 
0$ for $\nu(\varepsilon)\in (0,+\infty) \rightarrow 0$. 

The set $\displaystyle \Gamma_{1} = \nu^{-1}(0,1] $ is a semigroup. We note $\displaystyle 
\bar{\Gamma}_{1}= \Gamma_{1} \cup \left\{ 0\right\}$
On the set $\displaystyle 
\bar{\Gamma}= \Gamma \cup \left\{ 0\right\}$ we extend the operation on $\Gamma$ by adding the rules  
$00=0$ and $\varepsilon 0 = 0$ for any $\varepsilon \in \Gamma$. This is in agreement with the invariance 
of the end $0$ with respect to translations in $\Gamma$. 

In the Axioms A0, A1 we therefore may replace $\displaystyle [0,1]$ by $\displaystyle 
\bar{\Gamma}_{1}$, and so forth. 
\item[(b)] 
we may  leave some flexibility in Axiom A1 for the choice of base point of the dilatation, in the sense that 
\[\lim_{\nu(\varepsilon) \rightarrow 0}  \frac{1}{\nu(\varepsilon)}d(x , \delta^{x}_{\varepsilon} x) = 0\]
uniformly with respect to $x \in K$ compact set,  
 \item[(c)]
we may  relax the semigroup condition in the Axiom A2, in the sense: for any compact set $K \subset X$, for any  $x, \in K$, $\varepsilon, \mu$ with $\nu(\varepsilon), \nu(\mu)  \in (0,1)$ and $\displaystyle u, v \in \bar{B}_{d}(x,A)$   we have
\[ \frac{1}{\nu(\varepsilon \mu) }\mid d(\delta_{\varepsilon}^{x} \delta_{\mu}^{x} u , \delta_{\varepsilon}^{x} \delta_{\mu}^{x} v) - d(\delta_{\varepsilon \mu}^{x} u , \delta_{\varepsilon \mu}^{x} v) \mid \ \leq \  \mathcal{O}(\varepsilon \mu) \]
\item[(d)] 
in the Axioms A3 and A4 we may replace "$\varepsilon \rightarrow 0$" by "$\nu(\varepsilon) \rightarrow 0$" and "$1/\varepsilon$" by "$1/\nu(\varepsilon)$". 
\end{enumerate}

We shall write the proofs of further results   such that these work even if we modify the axioms in the 
sense explained above.  We shall nevertheless   use $\varepsilon$ and not $\nu(\varepsilon)$, 
in order to avoid  a too heavy notation. 

The axioms, as given in this section, are said to be in strong form. With the modifications explained at 
points (a), (b),  (c), (d) above, the axioms are said to be in weak form. 

Further, axioms are taken in weak form with the notational conventions explained above, unless it is explicitely 
stated that some axiom has to be taken in strong form. 

\subsection{Dilatation structures, tangent cones and metric profiles}

We shall explain now what the axioms mean. The first Axiom A1  is stating that the distance 
between $\displaystyle \delta^{x}_{\varepsilon} x$ and $x$ is negligible with respect to $\varepsilon$. 
If $\displaystyle \delta^{x}_{\varepsilon}x  = x$ then this axiom is trivially satisfied. 

The second Axiom A2. states that in an approximate sense the transformations $\displaystyle \delta^{x}_{\varepsilon}$ form an action of $\Gamma$ on $X$. As previously, if we suppose that 
\[\delta_{\varepsilon}^{x} \delta^{x}_{\mu} = \delta_{\varepsilon \mu}^{x}\] then this axiom is trivially 
satisfied. 

Remark now that  the binary tree formalism described in section 4 underlies and simplifies 
 the calculus with dilatation structures. 
More precisely, we shall use the results in section 4 in the proof of theorems in the next section.

The notation with binary trees for composition of dilatations is not directly adapted for taking limits as $\varepsilon 
\rightarrow 0$. An extension of the formalism can be made in this direction, but this would add length to this paper, 
which is devoted to first properties of dilatation structures. We reserve the full description of the formalism for a 
future paper.

In Axiom  A3 we take limits. In this subsection we shall look at dilatation structures from the metric point of view, 
by  using Gromov-Hausdorff distance and metric profiles. 
 
 We state the interpretation of the Axiom A3  as a theorem. 
But before a definition: we denote by $(\delta, \varepsilon)$ the distance on 
\[\bar{B}_{d^{x}}(x,1) = \left\{ y \in X \mbox{: } d^{x}(x,y) \leq 1 \right\}\] given by
\[(\delta, \varepsilon)(u,v) = \frac{1}{\varepsilon} d(\delta^{x}_{\varepsilon} u , \delta^{x}_{\varepsilon} v) \]

\begin{theorem}
\label{thcone}
Let $(X,d,\delta)$ be a dilatation structure. The following are consequences of the Axioms A0~-~A3 only: 
\begin{enumerate}
\item[(a)] 
for all $u,v \in X$ such that $\displaystyle d(x,u)\leq 1$ and $\displaystyle d(x,v) \leq 1$  and all $\mu \in (0,A)$ we have
\[d^{x}(u,v) = \frac{1}{\mu} d^{x}(\delta_{\mu}^{x} u , \delta^{x}_{\mu} v) \]
We shall say that $d^{x}$ has the cone property with respect to dilatations. 
\item[(b)] 
The curve $\displaystyle \varepsilon> 0 \mapsto \mathbb{P}^{x}(\varepsilon) = [\bar{B}_{d^{x}}(x,1), (\delta, \varepsilon), x]$ is a  metric profile.
\end{enumerate}
\end{theorem}

\begin{proof}
(a) For $\varepsilon, \mu \in (0,1)$ we have
\begin{align*}
\left| \frac{1}{\varepsilon \mu} d(\delta_{\varepsilon}^{x}\delta^{x}_{\mu} u, \delta_{\varepsilon}^{x}\delta^{x}_{\mu} v) - d^{x}(u,v) \right|
& \leq \,
\left| \frac{1}{\varepsilon \mu} d(\delta_{\varepsilon \mu}^{x} u, \delta_{\varepsilon}^{x}\delta^{x}_{\mu} u) -  
\frac{1}{\varepsilon \mu} d(\delta_{\varepsilon \mu}^{x} v, \delta_{\varepsilon}^{x}\delta^{x}_{\mu} v) \right|\\
& + \left|\frac{1}{\varepsilon \mu} d(\delta_{\varepsilon \mu}^{x}u, \delta_{\varepsilon \mu}^{x} v) - d^{x}(u,v)\right|
\end{align*}
Use now the Axioms A2 and A3 and pass to the limit with $\varepsilon \rightarrow 0$. This gives the desired equality. 

\medskip
(b)We have to prove that $\mathbb{P}^{x}$ is a metric profile. For this we have to compare two pointed metric spaces: 
\[ \left(\bar{B}_{d^{x}}(x,1), (\delta^{x}, \varepsilon \mu), x\right) \ \mbox{ and } \ \left(\bar{B}_{\frac{1}{\mu}(\delta^{x}, \varepsilon)}(x,1), \frac{1}{\mu}(\delta^{x}, \varepsilon), x \right) \]
Let $u \in X$ such that 
\[\frac{1}{\mu}(\delta^{x}, \varepsilon)(x, u) \leq 1 \]
This means that
\[\frac{1}{\varepsilon} d(\delta^{x}_{\varepsilon} x , \delta^{x}_{\varepsilon} u) \ \leq \ \mu \]
Further use the  Axioms  A1, A2 and the cone property proved before:
\[\frac{1}{\varepsilon} d^{x}(\delta^{x}_{\varepsilon} x, \delta^{x}_{\varepsilon} u) \ \leq \ (\mathcal{O}(\varepsilon) + 1) \mu\]
therefore, 
\[d^{x}(x,u)  \  \leq \ (\mathcal{O}(\varepsilon) + 1) \mu \]
It follows that for any $\displaystyle u \in \bar{B}_{\frac{1}{\mu}(\delta^{x}, \varepsilon)}(x,1)$ we can choose $\displaystyle w(u) \in \bar{B}_{d^{x}}(x,1)$ such that
\[\frac{1}{\mu} d^{x}(u, \delta^{x}_{\mu} w(u)) = \mathcal{O}(\varepsilon) \]
We want to prove that 
\[\mid \frac{1}{\mu}(\delta^{x}, \varepsilon) (u_{1}, u_{2}) - (\delta^{x}, \varepsilon \mu) (w(u_{1}), 
w(u_{2}) ) \mid \ \leq \ \mathcal{O}(\varepsilon \mu) + \frac{1}{\mu} \mathcal{O}(\varepsilon) + \mathcal{O}(\varepsilon) \ \]
This goes as follows: 
\begin{align*}
\lvert \frac{1}{\mu}(\delta^{x}, \varepsilon) (u_{1}, u_{2})  
&- (\delta^{x}, \varepsilon \mu) (w(u_{1}), w(u_{2}) ) \rvert
\,=\,
\left| \frac{1}{\varepsilon \mu} d( \delta^{x}_{\varepsilon} u_{1}, \delta^{x}_{\varepsilon} u_{2}) -  \frac{1}{\varepsilon \mu} d(\delta^{x}_{\varepsilon \mu} w(u_{1}), \delta^{x}_{\varepsilon \mu} w(u_{2}))\right|  \\
&\leq 
\mathcal{O}(\varepsilon \mu) + \left|\frac{1}{\varepsilon \mu} d( \delta^{x}_{\varepsilon} u_{1}, \delta^{x}_{\varepsilon} u_{2}) - \frac{1}{\varepsilon \mu} d(\delta^{x}_{\varepsilon}\delta^{x}_{ \mu} w(u_{1}), \delta^{x}_{\varepsilon}\delta^{x}_{ \mu} w(u_{2}))\right| \\
&\leq 
\mathcal{O}(\varepsilon \mu) + \frac{1}{\mu} \mathcal{O}(\varepsilon) \ + \ 
\frac{1}{\mu} \mid d^{x} (u_{1}, u_{2}) - d^{x}(\delta^{x}_{\mu} w(u_{1}) , \delta^{x}_{\mu} w(u_{2}))\mid
\end{align*}
In order to obtain the last estimate we used twice  the Axiom A3. We proceed as follows:  
\begin{gather*}
\mathcal{O}(\varepsilon \mu) \ + \ \frac{1}{\mu} \mathcal{O}(\varepsilon) \ + \ 
\frac{1}{\mu} \mid d^{x} (u_{1}, u_{2}) - d^{x}(\delta^{x}_{\mu} w(u_{1}) , \delta^{x}_{\mu} w(u_{2})) \mid \, \leq \ \\ 
\leq \mathcal{O}(\varepsilon \mu) \ + \ \frac{1}{\mu} \mathcal{O}(\varepsilon) \ + \ \frac{1}{\mu} 
d^{x}(u_{1}, \delta^{x}_{\mu} w(u_{1})) \  + \ \frac{1}{\mu} 
d^{x}(u_{1}, \delta^{x}_{\mu} w(u_{2}))\\ 
\leq \  \mathcal{O}(\varepsilon \mu) + \frac{1}{\mu} \mathcal{O}(\varepsilon) + \mathcal{O}(\varepsilon) 
\end{gather*}
This shows that the property (b) of a metric profile is satisfied. 
The property (a)  is proved in  the Theorem \ref{tmit}. 
\end{proof}

The following theorem is related to  Mitchell \cite{mitch} Theorem 1, concerning  sub-riemannian geometry. 

\begin{theorem}
\label{tmit}
In the hypothesis of theorem \ref{thcone}, we have  the following limit: 
\[\lim_{\varepsilon \rightarrow 0} \ \frac{1}{\varepsilon} \sup \left\{  \mid d(u,v) - d^{x}(u,v) \mid:  d(x,u) \leq \varepsilon, \ d(x,v) \leq \varepsilon \right\} = 0 \]
Therefore if $d^{x}$ is a true (i.e. nondegenerate) distance, then  $(X,d)$ admits a metric tangent space 
in $x$. 

Moreover, the metric profile $\displaystyle [\bar{B}_{d^{x}}(x,1), (\delta, \varepsilon), x]$ is  
almost nice, in the following sense. Let $c \in (0,1)$. Then  we have the inclusion
\[\delta^{x}_{\mu^{-1}}\left(\bar{B}_{\frac{1}{\mu}(\delta^{x}, \varepsilon)}(x,c)\right) \ \subset \ 
\bar{B}_{d^{x}}(x,1) \]
Moreover,  the following  Gromov-Hausdorff distance  is of order $\displaystyle \mathcal{O}(\varepsilon)$ 
for $\mu$ fixed (that is the modulus of convergence $\mathcal{O}(\varepsilon)$ does not depend on $\mu$): 
\[ \mu \ d_{GH}\left(  [\bar{B}_{d^{x}}(x,1) , (\delta^{x}, \varepsilon), x], \ [\delta^{x}_{\mu^{-1}}\left(\bar{B}_{\frac{1}{\mu}(\delta^{x}, \varepsilon)}(x,c)\right) , (\delta^{x}, \varepsilon \mu), x]\right) = 
\mathcal{O}(\varepsilon) \] 
For another Gromov-Hausdorff distance we have the estimate
\[d_{GH}\left( [\bar{B}_{\frac{1}{\mu}(\delta^{x}, \varepsilon)}(x,c), \frac{1}{\mu} (\delta^{x}, \varepsilon), x] \ , \ [\delta^{x}_{\mu^{-1}}\left(\bar{B}_{\frac{1}{\mu}(\delta^{x}, \varepsilon)}(x,c)\right) , (\delta^{x}, \varepsilon \mu), x]\right) = \mathcal{O}(\varepsilon \mu)\]
when $\varepsilon \in (0,\varepsilon(c))$. 
\end{theorem}

\begin{proof}
We start from the Axioms A0, A3 and we use the cone property. By A0,  for $\varepsilon \in (0,1)$ and $\displaystyle u,v \in  \bar{B}_{d}(x,\varepsilon)$ there exist $\displaystyle U,V  \in \bar{B}_{d}(x,A)$ such that \[u = \delta^{x}_{\varepsilon} U , v = \delta^{x}_{\varepsilon} V . \]
By the cone property we have 
\[\frac{1}{\varepsilon} \mid d(u,v) - d^{x}(u,v) \mid \,=\, \left|\frac{1}{\varepsilon} d(\delta^{x}_{\varepsilon} U,\delta^{x}_{\varepsilon} V) - d^{x}(U,V) \right| \]
By A2 we have 
\[  \left| \frac{1}{\varepsilon} d(\delta^{x}_{\varepsilon} U,\delta^{x}_{\varepsilon} V) - d^{x}(U,V) \right|\, \leq\, \mathcal{O}(\varepsilon)\]
This proves the first part of the theorem.

For the second part of the theorem take any 
$\displaystyle u \in \bar{B}_{\frac{1}{\mu}(\delta^{x}, \varepsilon)}(x,c)$. 
Then we have 
\[d^{x}(x,u) \, \leq \,  c \mu + \mathcal{O}(\varepsilon) \]
Then there exists   $\varepsilon(c) > 0$ such that 
for any $\varepsilon \in (0,\varepsilon(c))$ and $u$ in the mentioned ball we have
\[d^{x}(x,u) \,\leq \, \mu\]
In this case we can take directly $\displaystyle w(u) \,= \, \delta^{x}_{\mu^{-1}} u$ and simplify 
the  string of inequalities from the proof of  Theorem 
\ref{thcone}, point (b), to get eventually the three points 
from the second part of the theorem. 
\end{proof}

\section{Tangent bundle  of a dilatation structure}

In this section we shall use the calculus with binary decorated trees introduced in section 4, for a
space endowed with a dilatation structure. 

\subsection{Main results}

\begin{theorem} 
\label{thizo}
Let $(X,d,\delta)$ be a dilatation structure. Then the "infinitesimal translations" 
\[\displaystyle L^{x}_{u}(v) =  \lim_{\varepsilon \rightarrow 0}  \Delta^{x}_{\varepsilon}(u,v)\] are $\displaystyle d^{x}$ isometries. 
\end{theorem}

\begin{proof} 
The first part of the conclusion of  Theorem \ref{tmit} can be written as follows: 
\begin{equation}
\label{estiminf}
\sup \left\{ \frac{1}{\varepsilon} \mid d(u, v) - d^{x}(u,v) \mid \mbox{: } d(x,u) \leq \frac{3}{2}\varepsilon, \ 
d(x,v) \leq \frac{3}{2}\varepsilon \right\} \rightarrow 0
\end{equation}
as $\varepsilon \rightarrow 0$. 

 For $\varepsilon > 0$ sufficiently small the points 
$\displaystyle x, \delta^{x}_{\varepsilon}u, \delta^{x}_{\varepsilon} v , \delta^{x}_{\varepsilon} w$ 
are close one to another. Precisely, we have 
\[d(\delta^{\varepsilon}_{x} u, \delta^{\varepsilon}_{x} v) = \varepsilon( d^{x}(u,v) + \mathcal{O}(\varepsilon)) \]
Therefore, if we choose $u,v,w$ such that $d^{x}(u,v)< 1$ and  $d^{x}(u,w)< 1$, 
then there is $\eta>0$ such that for all $\varepsilon\in (0,\eta)$ we have 
\[d(\delta^{\varepsilon}_{x} u, \delta^{\varepsilon}_{x} v) \leq \frac{3}{2}\varepsilon,
 \qquad  d(\delta^{\varepsilon}_{x} u, 
\delta^{\varepsilon}_{x} v) \leq \frac{3}{2}\varepsilon \]

We  apply the estimate \eqref{estiminf}  for the 
basepoint $\displaystyle \delta^{x}_{\varepsilon} u$  to get
\[\frac{1}{\varepsilon} \mid d(\delta^{x}_{\varepsilon} v , \delta^{x}_{\varepsilon} w) - 
d^{\delta^{x}_{\varepsilon} u}(\delta^{x}_{\varepsilon} v , \delta^{x}_{\varepsilon} w) \mid \rightarrow 0\]
when $\varepsilon \rightarrow 0$. This can be written, using the cone property of the distance 
$\displaystyle d^{\delta^{x}_{\varepsilon} u}$, like
\begin{equation}
\label{esti2}
\left| \frac{1}{\varepsilon} d(\delta^{x}_{\varepsilon} v , \delta^{x}_{\varepsilon} w)  -
d^{\delta^{x}_{\varepsilon} u}\left( \delta_{\varepsilon^{-1}}^{\delta^{x}_{\varepsilon} u} \delta^{x}_{\varepsilon} v ,  \delta_{\varepsilon^{-1}}^{\delta^{x}_{\varepsilon} u} \delta^{x}_{\varepsilon} w 
\right) \right| \rightarrow 0
\end{equation}
as $\varepsilon \rightarrow 0$. By the Axioms A1, A3, the function 
\[(x,u,v) \mapsto d^{x}(u,v)\]
is an uniform limit of continuous functions, therefore uniformly continuous on compact sets. We can 
pass to the limit in the left hand side of the estimate \eqref{esti2}, using 
this uniform continuity and Axioms A3, A4, to get the result. 
\end{proof}

Let us define, in agreement with definition \ref{defsumdif} (b)
\[ \Sigma^{x}_{\varepsilon}(u, v) = \delta_{\varepsilon^{-1}}^{x} \delta_{\varepsilon}^{\delta_{\varepsilon}^{x} u} v\]

\begin{cor} 
\label{sigmaizo}
If for any $x$ the distance $\displaystyle d^{x}$ is non degenerate then there exists $C>0$ such that
for any $x$ and $u$ with $d(x,u) \leq C$ there exists a $\displaystyle d^{x}$ 
isometry $\displaystyle \Sigma^{x}(u, \cdot)$ obtained as the limit
\[ \lim_{\varepsilon \rightarrow 0} \Sigma^{x}_{\varepsilon}(u,v) = \Sigma^{x}(u, v)\] 
uniformly with respect to $x, u, v$ in compact set. 
\end{cor}

\begin{proof} 
From Theorem \ref{thizo} we know that $\displaystyle \Delta^{x}(u,\cdot)$ is  
a $\displaystyle d^{x}$ isometry. If $\displaystyle d^{x}$ is non degenerate then  $\displaystyle \Delta^{x}(u,\cdot)$ is invertible. Let  $\displaystyle \Sigma^{x}(u,\cdot)$ be the inverse. 

From Proposition \ref{psumdif} we know that $\displaystyle \Sigma^{x}_{\varepsilon}(u, \cdot)$ is the inverse of $\displaystyle \Delta^{x}_{\varepsilon}(u, \cdot)$. Therefore 
\begin{align*}
d^{x}( \Sigma^{x}_{\varepsilon}(u, w), \Sigma^{x}(u, w) ) 
&= d^{x}( \Delta^{x}(u,  \Sigma^{x}_{\varepsilon}(u, w)) , w)\\
&= d^{x}( \Delta^{x}(u,  \Sigma^{x}_{\varepsilon}(u, w)) ,  \Delta^{x}_{\varepsilon}(u,  \Sigma^{x}_{\varepsilon}(u, w))
\end{align*}
From the uniformity of convergence in Theorem \ref{thizo} and the uniformity assumptions in axioms of dilatation structures, the conclusion follows.
\end{proof}

The next theorem is the generalization of Proposition \ref{here3.4}. It is the main result of this paper. 

\begin{theorem}
\label{tgene}
Let $(X,d,\delta)$ be a dilatation structure which satisfies the strong form of the Axiom A2. Then for any $x \in X$ 
$\displaystyle (U(x), \Sigma^{x}, \delta^{x})$ is a conical group. Moreover, left translations of this group are 
$\displaystyle d^{x}$ isometries. 
\end{theorem}

\begin{proof} 
We start by proving that $\displaystyle (U(x), \Sigma^{x})$ is a local uniform group. The uniformities are induced by the distance $d$. 

We shall use the general relations written in terms of binary decorated trees. According to relation (\ref{difinv}) in Proposition \ref{p33}, we can pass to the limit with $\varepsilon \rightarrow 0$ and 
define
\[inv^{x}(u) = \lim_{\varepsilon \rightarrow 0} \Delta^{x}_{\varepsilon}(u,x) = \Delta^{x}(u,x) \]
From relation (\ref{invinv}) we get (after passing to the limit with $\varepsilon \rightarrow 0$) 
\[inv^{x}(inv^{x}(u)) = u\]
We shall see that $\displaystyle inv^{x}(u)$ is the inverse of $u$. 
Relation (\ref{eqinv}) gives 
\begin{equation}
\label{po}
\Delta^{x}(u,v) = \Sigma^{x} (inv^{x}(u), v)
\end{equation}
therefore  relations (a), (b) from Proposition \ref{psumdif} give 
\begin{gather}
\Sigma^{x}(inv^{x}(u), \Sigma^{x}(u,v)) = v 
\label{pb}\\
\Sigma^{x}(u, \Sigma^{x}(u,v)) = v
\label{pa}
\end{gather}
Relation (\ref{eqasoc2}) from Proposition \ref{p34} gives 
\begin{equation}
\label{pc}
\Sigma^{x}(u, \Sigma^{x}(v,w)) = \Sigma^{x}(\Sigma^{x}(u,v), w) 
\end{equation}
which shows that $\displaystyle \Sigma^{x}$ is an associative operation. From (\ref{pa}), (\ref{pb}) we 
obtain that for any $u,v$ 
\begin{gather}
\Sigma^{x}(\Sigma^{x}(inv^{x}(u), u),v) = v 
\label{pd}\\
\Sigma^{x}(\Sigma^{x}(u, inv^{x}(u)), v) = v
\label{pe}
\end{gather}
Remark that for any $x$, $v$ and $\varepsilon \in (0,1)$ we have $\displaystyle 
\Sigma^{x}(x,v) =  v$. indeed, this means that
\[ \begin{array}{ccccc}
\begin{parsetree}
 (..
    ( .$\bullet$.  
                                      .x.  
                                             (.$\circ$. (.$\circ$. .x. .x.)  .v.)
    ))\end{parsetree} &  = & \begin{parsetree}
 (..
    ( .$\bullet$.  
                                      .x.  
                                             (.$\circ$. .x.  .v.)
    ))\end{parsetree} &  = &  v
\end{array} \]
    Therefore $x$ is a neutral element at left for the operation $\displaystyle \Sigma^{x}$. 
From the definition of $\displaystyle inv^{x}$, relation (\ref{po}) and the fact that $\displaystyle inv^{x}$ 
is equal to its inverse, we get that $x$ is an inverse at right too: for any $x$, $v$ we have
\[\Sigma^{x}(v,x) = v \]
Replace now $v$ by $x$ in relations (\ref{pd}), (\ref{pe}) and prove that indeed $\displaystyle inv^{x}(u)$ 
is  the inverse of $u$.

We still have to prove that $\displaystyle (U(x), \Sigma^{x})$ admits $\displaystyle \delta^{x}$ as dilatations.In this reasoning we need the Axiom A2 in strong form. 

Namely we have to prove that for any $\mu \in (0,1)$ we have 
\[\delta^{x}_{\mu} \Sigma^{x}(u,v) = \Sigma^{x}(\delta^{x}_{\mu} u, \delta^{x}_{\mu} v) \]
For this is sufficient to notice that 
\[\Delta^{x}_{\varepsilon} \left( \delta^{x}_{\mu} u,  \delta^{x}_{\mu} v \right) =  \delta^{ \delta^{x}_{\epsilon \mu} u}_{\mu} \Delta^{x}_{\varepsilon \mu} (u,v) \]
and pass to the limit as $\varepsilon \rightarrow 0$. Notice that here we used the fact that dilatations 
$\displaystyle \delta^{x}_{\varepsilon}$ and $\displaystyle \delta^{x}_{\mu}$ exactly commute (Axiom A2 in strong form).  

Finally,  left translations $\displaystyle L^{x}_{u}$ are $\displaystyle d^{x}$ isometries. Really, this  is a straightforward consequence of  Theorem \ref{thizo} and corollary \ref{sigmaizo}. 
\end{proof}

The conical group $\displaystyle (U(x), \Sigma^{x}, \delta^{x})$ can be regarded as the tangent space of $(X,\delta,d)$ at $x$ and 
denoted further by $\displaystyle T_{x}X$. 

\subsection{Algebraic interpretation}

In order to better understand the algebraic structure of the sum, difference, 
inverse operations induced by a dilatation structure, we  collect 
previous results  regarding the properties of these operations, into one place. 

\begin{theorem}
\label{opcollection}
Let $(X,d,\delta)$ be a weak dilatation structure. Then, for any $x \in X$, $\varepsilon \in \Gamma$, 
$\nu(\varepsilon) < 1$,     we have 
\begin{enumerate}
\itemsep-3pt
\item[(a)] For any $u\in U(x)$,  $\displaystyle \Sigma_{\varepsilon}^{x}(x,u) = u$.
\item[(b)] For any $u\in U(x)$ the functions $\displaystyle  \Sigma_{\varepsilon}^{x}(u,\cdot)$ and 
$\displaystyle  \Delta_{\varepsilon}^{x}(u, \cdot)$ are inverse one to another. 
\item[(c)] The inverse function is shifted involutive:  for any $u\in U(x)$, 
\[inv^{\delta_{\varepsilon}^{x} u}_{\varepsilon} \, inv^{x}_{\varepsilon} (u) = u \]
\item[(d)] The sum operation is shifted associative: for  any $u, v, w$ sufficiently close to $x$ we have 
\[\Sigma_{\varepsilon}^{x} \left( u, \Sigma_{\varepsilon}^{\delta^{x}_{\varepsilon} u} (v, w)\right) = 
\Sigma_{\varepsilon}^{x} ( \Sigma^{x}(u,v), w) \]
\item[(e)] The difference, inverse and sum operations are related by 
\[ \Delta_{\varepsilon}^{x}(u,v) = \Sigma_{\varepsilon}^{\delta_{\varepsilon}^{x} u}
 \left( inv_{\varepsilon}^{x}(u), v\right) \]
 for any $u,v$ sufficiently close to $x$. 
 \item[(f)] For any $u,v$ sufficiently close to $x$ and $\mu \in \Gamma$, 
$\nu(\mu) < 1$,  we have
\[\Delta^{x}_{\varepsilon} \left( \delta^{x}_{\mu} u,  \delta^{x}_{\mu} v \right) =  \delta^{ \delta^{x}_{\epsilon \mu} u}_{\mu} \Delta^{x}_{\varepsilon \mu} (u,v) \]
\end{enumerate}
\end{theorem}

\section{Dilatation structures and differentiability}

\subsection{Equivalent dilatation structures}

\begin{definition}
\label{dilequi}
Two dilatation structures $(X, \delta , d)$ and $(X, \overline{\delta} , \overline{d})$  are equivalent  if 
\begin{enumerate}
\itemsep-3pt
\item[(a)] 
the identity  map $id: (X, d) \rightarrow (X, \overline{d})$ is bilipschitz  and 
\item[(b)]  
for any $x \in X$ there are functions $P^{x}, Q^{x}$ (defined for $u \in X$ sufficiently close to $x$) such that  
\begin{gather}
\lim_{\varepsilon \rightarrow 0} \frac{1}{\varepsilon} \overline{d} \left( \delta^{x}_{\varepsilon} u ,  \overline{\delta}^{x}_{\varepsilon} Q^{x} (u) \right)  = 0 
\label{dequiva}\\
\lim_{\varepsilon \rightarrow 0} \frac{1}{\varepsilon} d \left( \overline{\delta}^{x}_{\varepsilon} u ,  
\delta^{x}_{\varepsilon} P^{x} (u) \right)  = 0
\label{dequivb}
\end{gather}
uniformly with respect to $x$, $u$ in compact sets. 
\end{enumerate}
\end{definition}

\begin{proposition}
\label{pdilequi}
Two dilatation structures $(X, \delta , d)$ and $(X, \overline{\delta} , \overline{d})$  are equivalent  if and 
only if 
\begin{enumerate}
\itemsep-3pt
\item[(a)] 
the identity  map $id: (X, d) \rightarrow (X, \overline{d})$ is bilipschitz and 
\item[(b)]  
for any $x \in X$ there are functions $P^{x}, Q^{x}$ (defined for $u \in X$ sufficiently close to $x$) such that  
\begin{gather}
\lim_{\varepsilon \rightarrow 0}  \left(\overline{\delta}^{x}_{\varepsilon}\right)^{-1}  \delta^{x}_{\varepsilon} (u) = Q^{x}(u) 
\label{dequivap}\\
\lim_{\varepsilon \rightarrow 0}  \left(\delta^{x}_{\varepsilon}\right)^{-1}  \overline{\delta}^{x}_{\varepsilon} (u) = P^{x}(u) 
\label{dequivbp}
\end{gather}
uniformly with respect to $x$, $u$ in compact sets. 
\end{enumerate}
\end{proposition}

\begin{proof} 
We make the notations
\[
Q^{x}_{\varepsilon}(u) = \left(\overline{\delta}^{x}_{\varepsilon}\right)^{-1}  \delta^{x}_{\varepsilon} (u),
\qquad
P^{x}_{\varepsilon}(u) = \left(\delta^{x}_{\varepsilon}\right)^{-1}  \overline{\delta}^{x}_{\varepsilon} (u)
\]
The relation (\ref{dequiva}) is equivalent to 
\[
\mathcal{O}Ê(\varepsilon) + \overline{d}^{x} \left( Q^{x}_{\varepsilon}(u) ,  Q^{x}(u) \right) \rightarrow 0,
\qquad
\mathcal{O}Ê(\varepsilon) + d^{x} \left( P^{x}_{\varepsilon}(u) ,  P^{x}(u) \right) \rightarrow 0
\]
as $\varepsilon \rightarrow 0$, uniformly with respect to $x$, $u$ in compact sets. 
The conclusion follows after passing $\varepsilon \rightarrow 0$. 
\end{proof}
 
 The next theorem shows a link between the tangent bundles of equivalent dilatation structures. 
 
 \begin{theorem} 
 Let $(X, \delta , d)$ and $(X, \overline{\delta} , \overline{d})$  be  equivalent dilatation structures. Suppose that for any $x \in X$ the distance $d^{x}$ is non degenerate. Then for any $x \in X$ and 
 any $u,v \in X$ sufficiently close to $x$ we have:
 \begin{equation}
 \overline{\Sigma}^{x}(u,v) = Q^{x} \left( \Sigma^{x} \left( P^{x}(u) , P^{x}(v) \right)\right) 
 \label{isoequiv}
 \end{equation}
 The two tangent bundles  are therefore isomorphic in a natural sense. 
 \label{tisoequiv}
 \end{theorem}
 
 \begin{proof} 
 We notice first that the hypothesis is  symmetric: if $\displaystyle d^{x}$ is non degenerate then 
 $\displaystyle \overline{d}^{x}$ is non degenerate too. Really, this is 
straightforward from definition  \ref{dilequi} (a) and Axiom A3 for the two dilatation structures. 
 
 For the proof of relation (\ref{isoequiv}) is enough to remark that for $\varepsilon > 0$ but sufficiently small we have 
 \begin{equation}
  \overline{\Sigma}^{x}_{\varepsilon}(u,v) = Q^{x}_{\varepsilon} \left( \Sigma^{x}_{\varepsilon} \left(  P^{x}_{\varepsilon}(v) , P^{\overline{\delta}^{x}_{\varepsilon} u}_{\varepsilon} (v)  \right)\right) 
 \label{eqmo}
 \end{equation}
 Really, with tree notation, let 
\[\begin{array}{ccc}
 \begin{parsetree}
 (..
    ( .$\overline{\circ}$.  .x. .y.))
\end{parsetree} \  = \  \overline{\delta}^{x}_{\varepsilon} y ,  
\qquad 
\begin{parsetree}
 (..
    ( .$\circ$.  .x. .y.))
\end{parsetree} \  = \  \delta^{x}_{\varepsilon} y
    \end{array} \]
The relation (\ref{eqmo}), written from right to left, is
\[\begin{array}{ccc}
  \begin{parsetree}
 (..
    ( .$\overline{\bullet}$.  .x. ( .$\circ$.  .x. ( .$\bullet$.  .x.  ( .$\circ$.  ( .$\circ$.  .x. ( .$\bullet$.  .x. ( .$\overline{\circ}$.  .x. .u.))  )  ( .$\bullet$.  ( 
    .$\overline{\circ}$.  .x. .u.) ( .$\overline{\circ}$.  ( .$\overline{\circ}$.  .x. .u.) .v.))  ) )  )  ))
\end{parsetree}  & =  & \begin{parsetree}
 (..
    ( .$\overline{\bullet}$.  
                                      .x.  
                                             (.$\overline{\circ}$. (.$\overline{\circ}$. .x. .u.)  .v.)
    ))\end{parsetree}
    \end{array} \]
But this is true by cancellations  of dilatations and definitions of the operators $\displaystyle P^{x}_{\varepsilon}$ and $\displaystyle Q^{x}_{\varepsilon}$. 
\end{proof}

\subsection{Differentiable functions}

Dilatation structures allow to define differentiable functions.
The idea is to keep only  one relation from definition \ref{dilequi}, 
namely (\ref{dequiva}). We also renounce to uniform convergence with respect 
 to $x$ and $u$, and we replace this with uniform convergence in $u$ and with 
a conical group morphism condition for the derivative. 

First we need the natural definition below. 

\begin{definition}
\label{defmorph}
 Let $(N,\delta)$ and $(M,\bar{\delta})$ be two conical groups. A function $f:N\rightarrow M$ 
is a conical group morphism if $f$ is a group morphism and for any $\varepsilon>0$ and $u\in N$ we have 
 $\displaystyle f(\delta_{\varepsilon} u) = \bar{\delta}_{\varepsilon} f(u)$. 
\end{definition}

The definition of derivative with respect to dilatations structures follows. 

 \begin{definition}
\label{defdif} 
Let $(X, \delta , d)$ and $(Y, \overline{\delta} , \overline{d})$ be two dilatation structures and $f:X \rightarrow Y$ be a continuous function. The function $f$ is differentiable in $x$ if there exists a 
 conical group morphism  $\displaystyle Q^{x}:T_{x}X\rightarrow T_{f(x)}Y$, defined on a neighbourhood of $x$ with values in  a neighbourhood  of $f(x)$ such that 
\begin{equation}
\label{edefdif}
\lim_{\varepsilon \rightarrow 0} \sup \left\{  \frac{1}{\varepsilon} \overline{d} \left( f\left( \delta^{x}_{\varepsilon} u\right) ,  \overline{\delta}^{f(x)}_{\varepsilon} Q^{x} (u) \right) \mbox{: } d(x,u) \leq \varepsilon \right\} = 0 
\end{equation}
The morphism $\displaystyle Q^{x}$ is called the derivative of $f$ at $x$ and will be sometimes denoted by $Df(x)$.

The function $f$ is uniformly differentiable if it is differentiable everywhere and the limit in (\ref{edefdif}) 
is uniform in $x$ in compact sets. 
\end{definition}
 
This definition deserves a short discussion. 
 Let $(X, \delta , d)$ and $(Y, \overline{\delta} , \overline{d})$ be two 
dilatation structures and $f:X \rightarrow Y$  a function differentiable in 
$x$. The derivative of $f$ in $x$ is a 
 conical group morphism  $\displaystyle Df(x):T_{x}X\rightarrow T_{f(x)}Y$, 
which means that $Df(x)$ is 
 defined on a open set around $x$ with values in a open set around $f(x)$, 
having the following properties:
 \begin{enumerate}
 \item[(a)] for any $u,v$ sufficiently close to $x$ 
 \[Df(x)\left(\Sigma^{x}(u,v)\right) = \Sigma^{f(x)}\left(Df(x)(u), Df(x)(v)\right)\]
 \item[(b)] for any $u$ sufficiently close to $x$ and any $\varepsilon \in (0,1]$ 
 \[Df(x)\left(\delta^{x}_{\varepsilon} u\right) = \bar{\delta}^{f(x)}_{\varepsilon}\left(Df(x)(u)\right) \]
 \item[(c)] the function $Df(x)$ is continuous, as uniform limit of 
continuous functions. Indeed, the relation 
 (\ref{edefdif}) is equivalent to the existence of the uniform limit (with respect to $u$ in compact sets)
 \[
Df(x)(u) = \lim_{\varepsilon\rightarrow 0} \bar{\delta}^{f(x)}_{\varepsilon^{-1}} \left( f\left( \delta_{\varepsilon}^{x} u\right)\right)
\] 
 \end{enumerate}
 
 From (\ref{edefdif}) alone and axioms of dilatation structures we can prove 
properties (b) and (c). 
 We can reformulate therefore the definition of the derivative by asking that 
$Df(x)$ exists as an uniform 
 limit (as in point (c) above) and that $Df(x)$ has the property (a) above. 
From these considerations the chain rule for derivatives is straightforward.

A trivial way to obtain a differentiable function (everywhere)  is to modify the dilatation 
structure on the target space. 
 
 \begin{definition}
\label{ddif} 
Let  $(X, \delta , d)$ be a  dilatation structure and $f:(X, d) \rightarrow (Y, \overline{d})$ be a bilipschitz  and surjective  function. We define then the transport of $(X, \delta , d)$ by $f$, named $(Y, f*\delta , \overline{d})$, by 
 \[\left( f*\delta\right)^{f(x)}_{\varepsilon} f(u) = f \left( \delta^{x}_{\varepsilon} u \right)\]
 \end{definition}
 
The relation of differentiability with equivalent dilatation structures is given by the following simple
\begin{proposition}
 Let  $(X, \delta , d)$ and $(X, \overline{\delta} , \overline{d})$ be two 
dilatation structures and $f:(X, d) \rightarrow (X, \overline{d})$ be a 
bilipschitz  and surjective  function. The dilatation structures 
$(X, \overline{\delta} , \overline{d})$ and $(X, f*\delta , \overline{d})$ are 
equivalent if and only if $f$ and $\displaystyle f^{-1}$ are uniformly  
differentiable. 
 \label{peqd}
 \end{proposition}

 \begin{proof} 
 Straightforward from  definitions \ref{dilequi}   and \ref{ddif}. 
\end{proof}

\section{Differential structure, conical groups and dilatation structures} 

In this section we collect some facts which relate differential structures with dilatation structures. 
We resume then the paper with a justification of the unusual way of defining uniform groups (definition 
 \ref{dunifg}) by the fact that the $op$ function (the group operation) is differentiable with respect to dilatation structures which are natural for a group with dilatations. 
 
 \subsection{Differential structures and dilatation structures}
  
A differential structure on a manifold is an equivalence class of compatible atlases. We show here that an atlas induces an 
equivalence class of  dilatation structures and that two compatible atlases induce the same equivalence class of  dilatation 
structures. 

Let $M$ be a $\displaystyle \mathcal{C}^{1}$ $n$-dimensional real manifold and  
$\displaystyle \mathcal{A}$ an atlas of this manifold. For each chart 
$\displaystyle \phi : W \subset M \rightarrow \mathbb{R}^{n}$  we shall define a dilatation structure on $W$. 

Suppose that $\displaystyle \phi(W) \subset \mathbb{R}^{n}$ is convex (if not then take an open subset of $W$ with 
this property).  For $x, u \in W$ and $\varepsilon \in (0,1]$ define the dilatation 
\[\delta^{x}_{\varepsilon} u = \phi^{-1} \left( \phi(x) + \varepsilon (\phi(u) -\phi(x)) \right) \]
Otherwise said, the dilatations in $W$ are transported from $\displaystyle \mathbb{R}^{n}$. Equally, we transport on $W$ the euclidean distance of $\displaystyle \mathbb{R}^{n}$. We obviously get a dilatation structure on $W$. 

If we have two charts 
$\displaystyle \phi_{i} : W_{i} \subset M \rightarrow \mathbb{R}^{n}$, $i=1,2$, 
belonging to the same atlas $\displaystyle \mathcal{A}$, then we have two  equivalent dilatation structures on 
$\displaystyle W_{1} \cap W_{2}$. Indeed, the atlas $\displaystyle \mathcal{A}$ is $\displaystyle \mathcal{C}^{1}$ therefore 
the distances (induced from the charts) are (locally) in bilipschitz equivalence. 
Denote by $\displaystyle  \overline{\delta}$  the dilatation obtained from the chart $\displaystyle \phi_{2}$.  
A short computation shows that (we use here the transition map 
$\displaystyle \phi_{21} = \phi_{2} (\phi_{1})^{-1}$)
\[Q^{x}_{\varepsilon}(u) = (\phi_{2})^{-1} \left( \phi_{2}(x) + \frac{1}{\varepsilon} \left( \phi_{21} \left( \phi_{1} (x) + \varepsilon ( f(u) - f(x) ) \right) - \phi_{2}(x) \right) \right)\]
 therefore, as $\varepsilon \rightarrow 0$, we have 
 \[\lim_{\varepsilon \rightarrow 0} Q^{x}_{\varepsilon} (u) = Q^{x}(u) = (\phi_{2})^{-1} \left( \phi_{2}(x) + 
 D \phi_{21}(f(x))( f(u) - f(x) )  \right)\]
 A similar computation shows that $\displaystyle P^{x}$ also exists. The uniform convergence requirements come from the fact that we use a $\displaystyle \mathcal{C}^{1}$ atlas. 
 
 A similar reasoning shows that in fact two compatible atlases induce the same equivalence class of  dilatation structures. 
 
 \subsection{Conical groups and dilatation structures}
 
 In a group with dilatations $(G, \delta)$  we define dilatations based in any point $x \in G$ by 
 \begin{equation}
 \delta^{x}_{\varepsilon} u = x \delta_{\varepsilon} ( x^{-1}u)
 \label{dilat}
 \end{equation}
 
 \begin{definition} 
\label{dnco}
A normed group with dilatations $(G, \delta, \| \cdot \|)$ is a 
group with dilatations  $(G, \delta)$ endowed with a continuous norm  
function $\displaystyle \|\cdot \|: G \rightarrow \mathbb{R}$ which satisfies 
(locally, in a neighbourhood  of the neutral element $e$) the following properties: 
 \begin{enumerate}
\itemsep-2pt
 \item[(a)] for any $x$ we have $\| x\| \geq 0$; if $\| x\| = 0$ then $x=e$, 
 \item[(b)] for any $x,y$ we have $\|xy\| \leq \|x\| + \|y\|$, 
 \item[(c)] for any $x$ we have $\displaystyle \| x^{-1}\| = \|x\|$, 
 \item[(d)] the limit 
$\displaystyle \lim_{\varepsilon \rightarrow 0} \frac{1}{\nu(\varepsilon)} \| \delta_{\varepsilon} x \| = \| x\|^{N}$ 
 exists, is uniform with respect to $x$ in compact set, 
 \item[(e)] if $\displaystyle \| x\|^{N} = 0$, then $x=e$.
\end{enumerate}
\end{definition}

It is easy to see that if $(G, \delta, \| \cdot \|)$ is a normed group with dilatations
then $\displaystyle (G, \beta, \delta, \|\cdot\|^{N})$ is a normed conical group. 
The norm $\displaystyle \|\cdot\|^{N}$ satisfies the 
stronger form of property (d) of Definition \ref{dnco}: for any $\varepsilon >0$,
$\| \delta_{\varepsilon} x \|^{N} = \varepsilon \| x \|^{N}$.

Normed groups with dilatations can be encountered in sub-Riemannian geometry. Normed conical groups generalize the notion of Carnot groups. 

In a normed group with dilatations we have a natural left invariant distance given by
\begin{equation}
d(x,y) = \| x^{-1}y\|
\label{dnormed}
\end{equation}

\begin{theorem}
\label{tgrd}
Let $(G, \delta, \| \cdot \|)$ be  a locally compact  normed group with dilatations. Then $(G, \delta, d)$ is 
a dilatation structure, where $\delta$ are the dilatations defined by (\ref{dilat}) and the distance $d$ is induced by the norm as in (\ref{dnormed}). 
\end{theorem}

\begin{proof}
The Axiom A0 is straightforward from definition \ref{dunifg}, 
definition \ref{defgwd},  Axiom H0,  and because the dilatation structure is left invariant, in the sense that the transport by left translations in $G$, according to Definition \ref{ddif}, preserves the dilatations $\delta$. We also trivially have Axioms A1 and A2 satisfied. 

For the Axiom A3  remark that
\[d(\delta_{\varepsilon}^{x} u , \delta_{\varepsilon}^{x} v) = d(x \delta_{\varepsilon}(x^{-1}u),  x \delta_{\varepsilon}(x^{-1}u)) = d( \delta_{\varepsilon}(x^{-1}u) ,  \delta_{\varepsilon}(x^{-1}v)) \]
Denote $\displaystyle U= x^{-1}u$, $\displaystyle V= x^{-1}v$ and for $\varepsilon > 0$ let 
\[\displaystyle \beta_{\varepsilon}(u, v) =  \delta_{\varepsilon}^{-1}
\left((\delta_{\varepsilon}u) (\delta_{\varepsilon}v ) \right)\]
 We have then: 
\[\frac{1}{\varepsilon} d(\delta_{\varepsilon}^{x} u , \delta_{\varepsilon}^{x} v) = 
\frac{1}{\varepsilon} \|\delta_{\varepsilon} \beta_{\varepsilon} \left(  \delta_{\varepsilon}^{-1}
\left( ( \delta_{\varepsilon}V)^{-1}\right) , U\right)\| \]

Define the function 
\[d^{x}(u,v) = \| \beta( V^{-1} , U) \|^{N} \]
From Definition \ref{defgwd} Axioms H1, H2, and from definition \ref{dnco} (d), we obtain that Axiom A3 
is satisfied. 

For the Axiom A4 we have to calculate
\begin{align*}
\Delta^{x}(u,v) 
&= \delta_{\varepsilon^{-1}}^{\delta^{x}_{\varepsilon} u} \delta^{x}_{\varepsilon} v \\
&= \left(\delta^{x}_{\varepsilon} u \right) \left( \delta_{\varepsilon}\right)^{-1} \left( \left(\delta^{x}_{\varepsilon} u \right)^{-1} \left(\delta^{x}_{\varepsilon} v \right)\right)  \\
&= 
\left(x\delta_{\varepsilon} U \right)  \beta_{\varepsilon}\left(  \delta_{\varepsilon}^{-1}
\left( ( \delta_{\varepsilon}V)^{-1}\right) , U\right) \rightarrow x \beta\left( V^{-1}, U \right)  
\end{align*}
as $\varepsilon \rightarrow 0$. Therefore the Axiom A4 is satisfied.
\end{proof}

We remarked in the proof of the previous theorem that the transport by left translations in $G$, according to 
Definition \ref{ddif}, preserves the dilatation structure on $G$. This implies, according to Proposition 
\ref{peqd}, that left translations are differentiable. On the contrary, a short computation and examples 
from sub-Riemannian geometry indicate that right translations are not differentiable. 

Nevertheless, the operation $op$ is differentiable, if we endow the group $\displaystyle G^{(2)} = G \times G$ with a good dilatation structure. This will justify the non standard way to define local uniform 
groups in Definition \ref{dunifg}. 

Start from the fact that if $G$ is a local uniform group then $\displaystyle G^{(2)} $ is a local uniform group too. If $G$ is also normed, with dilatations, then we can easily define a similar structure on 
$\displaystyle G^{(2)}$. Really, the norm on $\displaystyle G^{(2)}$ can be taken as 
\[\|(x,y)\|^{(2)} = \max \left\{ \|x\| , \|y\| \right\}\]
and dilatations 
\[\delta_{\varepsilon}^{(2)} (x,y) = (\delta_{\varepsilon} x , \delta_{\varepsilon} y )\]
We leave to the reader to check that $\displaystyle G^{(2)}$ endowed with this norm and these dilatations is indeed a normed group with dilatations. 

\begin{theorem}
\label{opsm}
Let $(G, \delta, \| \cdot \|)$ be  a locally compact  normed group with dilatations and let
$(G^{(2)}, \delta^{(2)}, \|\cdot\|^{(2)})$ be the associated normed group with dilatation. 
Then  the operation ($op$ function) is differentiable. 
\end{theorem}

\begin{proof}
We start from the formula (easy to check in $\displaystyle G^{(2)}$)
\[(x,y)^{-1} = (x^{-1}, x y^{-1} x^{-1}) \]
Then we have 
\[\delta^{(x,y)}_{\varepsilon} (u,v) =  \left(x \delta_{\varepsilon} (x^{-1} u) , \left(  
\delta_{\varepsilon} (x^{-1} u)\right)^{-1} y   \delta_{\varepsilon} (x^{-1} u) 
\delta_{\varepsilon} \left( u^{-1} xy^{-1} x^{-1} uv\right) \right) \]
Let us define 
\[Q^{(x,y)}(u,v) = op(x,y) \beta((x,y)^{-1}(u,v))\]
Then we have 
\[
\frac{1}{\varepsilon} d\left( op\left( \delta^{(x,y)}(u,v) \right) , \delta^{op(x,y)} Q^{(x,y)}(u,v) \right) 
= 
\frac{1}{\varepsilon} d \left( \delta_{\varepsilon} \beta_{\varepsilon} ((x,y)^{-1}(u,v)) , \delta_{\varepsilon} \beta((x,y)^{-1}(u,v)) \right) 
\]
The right hand side of this equality  converges then to $0$ as $\varepsilon \rightarrow 0$. More precisely, we have 
\begin{gather*}
\sup \left\{  \frac{1}{\varepsilon} d\left( op\left( \delta^{(x,y)}(u,v) \right) , \delta^{op(x,y)} Q^{(x,y)}(u,v) \right)  \mbox{: }Êd^{(2)}((x,y), (u,v)) \leq \varepsilon \right\} =\\
= 
\sup \left\{ d^{e} \left( \beta_{\varepsilon} ((x,y)^{-1}(u,v)) , \beta ((x,y)^{-1}(u,v)) \right) \mbox{: } 
d^{(2)}((x,y), (u,v)) \leq \varepsilon \right\} + \mathcal{O}(\varepsilon) 
 \tag*{\qed}
\end{gather*}
 \renewcommand{\qed}{}
\end{proof}

In particular, we have $\displaystyle Q^{(e,e)}(u,v) = \beta(u,v)$, which shows that the operation $\beta$
is the differential of the operation $op$ computed in the neutral element of $\displaystyle G^{(2)}$.

\smallskip
\Date{Received October 14, 2006\\Revised February 25, 2007}

\label{lastpage}

\end{document}